# HIGH-DIMENSIONAL CLASSIFICATION USING FEATURES ANNEALED INDEPENDENCE RULES[1]

By Jianqing Fan and Yingying Fan

*Princeton University and Harvard University*

Classification using high-dimensional features arises frequently in many contemporary statistical studies such as tumor classification using microarray or other high-throughput data. The impact of dimensionality on classifications is poorly understood. In a seminal paper, Bickel and Levina [*Bernoulli* **10** (2004) 989–1010] show that the Fisher discriminant performs poorly due to diverging spectra and they propose to use the independence rule to overcome the problem. We first demonstrate that even for the independence classification rule, classification using all the features can be as poor as the random guessing due to noise accumulation in estimating population centroids in high-dimensional feature space. In fact, we demonstrate further that almost all linear discriminants can perform as poorly as the random guessing. Thus, it is important to select a subset of important features for high-dimensional classification, resulting in Features Annealed Independence Rules (FAIR). The conditions under which all the important features can be selected by the two-sample *t*-statistic are established. The choice of the optimal number of features, or equivalently, the threshold value of the test statistics are proposed based on an upper bound of the classification error. Simulation studies and real data analysis support our theoretical results and demonstrate convincingly the advantage of our new classification procedure.

**1. Introduction.** With rapid advance of imaging technology, high-throughput data such as microarray and proteomics data are frequently seen in many contemporary statistical studies. For instance, in the analysis of Microarray data, the dimensionality is frequently thousands or more, while the sample

---

Received January 2007; revised May 2007.

[1]Supported by NSF Grants DMS-03-54223 and DMS-07-04337 and NIH Grant R01-GM072611.

size is typically in the order of tens [West et al. (2001) and Dudoit, Fridlyand and Speed (2002)]; see Fan and Ren (2006) for an overview. The large number of features presents an intrinsic challenge to classification problems. For an overview of statistical challenges associated with high dimensionality, see Fan and Li (2006).

Classical methods of classification break down when the dimensionality is extremely large. For example, even when the covariance matrix is known, Bickel and Levina (2004) demonstrate convincingly that the Fisher discriminant analysis performs poorly in a minimax sense due to the diverging spectra (e.g., the condition number goes to infinity as dimensionality diverges) frequently encountered in the high-dimensional covariance matrices. Even if the true covariance matrix is not ill conditioned, the singularity of the sample covariance matrix will make the Fisher discrimination rule inapplicable when the dimensionality is larger than sample size. Bickel and Levina (2004) show that the independence rule overcomes the above two problems. However, in tumor classification using microarray data, we hope to find tens of genes that have high discriminative power. The independence rule, studied by Bickel and Levina (2004), does not possess this kind of properties.

The difficulty of high-dimensional classification is intrinsically caused by the existence of many noise features that do not contribute to the reduction of misclassification rate. Though the importance of dimension reduction and feature selection has been stressed and many methods have been proposed in the literature, very little research has been done on theoretical analysis of the impacts of high dimensionality on classification. For example, using most discrimination rules such as the linear discriminants, we need to estimate the population mean vectors from the sample. When the dimensionality is high, even though each component of the population mean vectors can be estimated with accuracy, the aggregated estimation error can be very large and this has adverse effects on the misclassification rate. Therefore, when there is only a fraction of features that account for most of the variation in the data such as tumor classification using gene expression data, using all features will increase the misclassification rate.

To illustrate the idea, we study independence classification rule. Specifically, we give an explicit formula on how the signal and noise affect the misclassification rates. We show formally how large the signal to noise ratio can be such that the effect of noise accumulation can be ignored, and how small this ratio can be before the independence classifier performs as poorly as the random guessing. Indeed, as demonstrated in Section 2, the impact of the dimensionality can be very drastic. For the independence rule, the misclassification rate can be as high as the random guessing even when the problem is perfectly classifiable. In fact, we demonstrate that almost all linear discriminants cannot perform any better than random guessing, due to the noise accumulation in the estimation of the population mean vectors, unless the signals are very strong, namely the population mean vectors are very far apart.

The above discussion reveals that feature selection is necessary for high-dimensional classification problems. When the independence rule is applied to selected features, the resulting Feature Annealed Independent Rules (FAIR) overcome both the issues of interpretability and the noise accumulation. One can extract the important features via variable selection techniques such as the penalized quasi-likelihood function. See Fan and Li (2006) for an overview. One can also employ a simple two-sample $t$-test as in Tibshirani et al. (2002) to identify important genes for the tumor classification, resulting in the nearest shrunken centroids method. Such a simple method corresponds to a componentwise regression method or a ridge regression method with ridge parameters tending to $\infty$ [Fan and Lv (2007)]. Hence, it is a specific and useful example of the penalized quasi-likelihood method for feature selection. It is surprising that such a simple proposal can indeed extract all important features. Indeed, we demonstrate that under suitable conditions, the two-sample $t$-statistic can identify all the features that efficiently characterize both classes.

Another popular class of the dimension reduction methods is projection. They have been widely applied to the classification based on the gene expression data. See, for example, principal component analysis in Ghosh (2002), Zou, Hastie and Tibshirani (2004) and Bair et al. (2006); partial least squares in Nguyen and Rocke (2002), Huang and Pan (2003) and Boulesteix (2004); and sliced inverse regression in Chiaromonte and Martinelli (2002), Antoniadis, Lambert-Lacroix and Leblanc (2003) and Bura and Pfeiffer (2003). These projection methods attempt to find directions that can result in small classification errors. In fact, the directions found by these methods usually put much more weight on features that have large classification power. In general, however, linear projection methods are likely to perform poorly unless the projection vector is sparse, namely, the effective number of selected features is small. This is due to the aforementioned noise accumulation prominently featured in high-dimensional problems, recalling discrimination based on linear projections onto almost all directions can perform as poorly as the random guessing.

As direct application of the independence rule is not efficient, we propose a specific form of FAIR. Our FAIR selects the statistically most significant $m$ features according to the componentwise two-sample $t$-statistics between two classes, and applies the independence classifiers to these $m$ features. Interesting questions include how to choose the optimal $m$, or equivalently, the threshold value of $t$-statistic, such that the classification error is minimized, and how this classifier performs compared with the independence rule without feature selection and the oracle-assisted FAIR. All these questions will be formally answered in this paper. Surprisingly, these results are similar to those for the adaptive Neyman test in Fan (1996). The theoretical results also indicate that FAIR without oracle information performs worse than the one with oracle information, and the difference of classification error depends on the threshold value, which is consistent with the common sense.

There is a huge literature on classification. To name a few in addition to those mention before, Bai and Saranadasa (1996) dealt with the effect of high dimensionality in a two-sample problem from a hypothesis testing viewpoint; Friedman (1989) proposed a regularized discriminant analysis to deal with the problems associated with high dimension while performing computations in the regular way; Dettling and Bühlmann (2003) and Bühlmann and Yu (2003) study boosting with logit loss and $L_2$ loss, respectively, and demonstrate the good performances of these methods in high-dimensional setting; Greenshtein and Ritov (2004), Greenshtein (2006) and Meinshausen (2007) introduced and studied the concept of persistence, which places more emphasis on misclassification rates or expected loss rather than the accuracy of estimated parameters.

This article is organized as follows. In Section 2, we demonstrate the impact of dimensionality on the independence classification rule, and show that discrimination based on projecting observations onto almost all linear directions is nearly the same as random guessing. We establish, in Section 3, the conditions under which two-sample $t$-test can identify all the important features with probability tending to 1. In Section 4, we propose FAIR and give an upper bound of its classification error. Simulation studies and real data analyses are conducted in Section 5. The conclusion of our study is summarized in Section 6. All proofs are given in the Appendix.

**2. Impact of high dimensionality.**  Consider the $p$-dimensional classification problem between two classes $\mathcal{C}_1$ and $\mathcal{C}_2$. Suppose that from class $\mathcal{C}_k$, we have $n_k$ observations $\mathbf{Y}_{k1}, \ldots, \mathbf{Y}_{kn_k}$ in $\mathbb{R}^p$. The $j$th feature of the $i$th sample from class $\mathcal{C}_k$ satisfies the model

$$(2.1) \qquad Y_{kij} = \mu_{kj} + \epsilon_{kij}, \qquad k = 1, 2, \ i = 1, \ldots, n_k, \ j = 1, \ldots, p,$$

where $\mu_{kj}$ is the mean effect of the $j$th feature in class $\mathcal{C}_k$ and $\epsilon_{kij}$ is the corresponding Gaussian random noise for $i$th observation. In matrix notation, the above model can be written as

$$\mathbf{Y}_{ki} = \boldsymbol{\mu}_k + \boldsymbol{\epsilon}_{ki}, \qquad k = 1, 2, \ i = 1, \ldots, n_k,$$

where $\boldsymbol{\mu}_k = (\mu_{k1}, \ldots, \mu_{kp})'$ is the mean vector of class $\mathcal{C}_k$ and $\boldsymbol{\epsilon}_{ki} = (\epsilon_{ki1}, \ldots, \epsilon_{kip})'$ has the distribution $N(\mathbf{0}, \boldsymbol{\Sigma}_k)$. We assume that all observations are independent across samples and in addition, within class $\mathcal{C}_k$, observations $\mathbf{Y}_{k1}, \ldots, \mathbf{Y}_{kn_k}$ are also identically distributed. Throughout this paper, we make the assumption that the two classes have compatible sample sizes, that is, $c_1 \le n_1/n_2 \le c_2$ with $c_1$ and $c_2$ some positive constants.

We first investigate the impact of high dimensionality on classification. For simplicity, we temporarily assume that the two classes $\mathcal{C}_1$ and $\mathcal{C}_2$ have the same covariance matrix $\boldsymbol{\Sigma}$. To illustrate our idea, we consider the independence classification rule, which classifies the new feature vector $\mathbf{x}$ into class $\mathcal{C}_1$ if

$$\delta(\mathbf{x}) = (\mathbf{x} - \boldsymbol{\mu})' \mathbf{D}^{-1} (\boldsymbol{\mu}_1 - \boldsymbol{\mu}_2) > 0,$$

where $\boldsymbol{\mu} = (\boldsymbol{\mu}_1 + \boldsymbol{\mu}_2)/2$ and $\mathbf{D} = \operatorname{diag}(\boldsymbol{\Sigma})$. This classifier has been thoroughly studied in Bickel and Levina (2004). They showed that in the classification of two normal populations, this independence rule greatly outperforms the Fisher linear discriminant rule under broad conditions when the number of variables is large.

The independence rule depends on the marginal parameters $\boldsymbol{\mu}_1$, $\boldsymbol{\mu}_2$ and $\mathbf{D} = \operatorname{diag}\{\sigma_1^2, \ldots, \sigma_p^2\}$. They can easily be estimated from the samples

$$\widehat{\boldsymbol{\mu}}_k = \sum_{i=1}^{n_k} \mathbf{Y}_{ki}/n_k, \qquad k = 1, 2, \qquad \widehat{\boldsymbol{\mu}} = (\widehat{\boldsymbol{\mu}}_1 + \widehat{\boldsymbol{\mu}}_2)/2$$

and

$$\widehat{\mathbf{D}} = \operatorname{diag}\{(S_{1j}^2 + S_{2j}^2)/2, \ j = 1, \ldots, p\},$$

where $S_{kj}^2 = \sum_{i=1}^{n_k} (Y_{kij} - \bar{Y}_{kj})^2/(n_k - 1)$ is the sample variance of the $j$th feature in class $k$ and $\bar{Y}_{kj} = \sum_{i=1}^{n_k} Y_{ki}/n_k$. Hence, the plug-in discrimination function is

$$\hat{\delta}(\mathbf{x}) = (\mathbf{x} - \widehat{\boldsymbol{\mu}})' \widehat{\mathbf{D}}^{-1} (\widehat{\boldsymbol{\mu}}_1 - \widehat{\boldsymbol{\mu}}_2).$$

Denote the parameter by $\boldsymbol{\theta} = (\boldsymbol{\mu}_1, \boldsymbol{\mu}_2, \boldsymbol{\Sigma})$. If we have a new observation $\mathbf{X}$ from class $\mathcal{C}_1$, then the misclassification rate of $\hat{\delta}$ is

$$(2.2) \quad W(\hat{\delta}, \boldsymbol{\theta}) = P(\hat{\delta}(\mathbf{X}) \le 0 | \mathbf{Y}_{ki}, i = 1, \ldots, n_k, \ k = 1, 2) = 1 - \Phi(\Psi),$$

where

$$\Psi = \frac{(\boldsymbol{\mu}_1 - \widehat{\boldsymbol{\mu}})' \widehat{\mathbf{D}}^{-1} (\widehat{\boldsymbol{\mu}}_1 - \widehat{\boldsymbol{\mu}}_2)}{\sqrt{(\widehat{\boldsymbol{\mu}}_1 - \widehat{\boldsymbol{\mu}}_2)' \widehat{\mathbf{D}}^{-1} \boldsymbol{\Sigma} \widehat{\mathbf{D}}^{-1} (\widehat{\boldsymbol{\mu}}_1 - \widehat{\boldsymbol{\mu}}_2)}},$$

and $\Phi(\cdot)$ is the standard Gaussian distribution function. The worst case classification error is

$$W(\hat{\delta}) = \max_{\boldsymbol{\theta} \in \Gamma} W(\hat{\delta}, \boldsymbol{\theta}),$$

where $\Gamma$ is some parameter space to be defined. Let $n = n_1 + n_2$. In our asymptotic analysis, we always consider the misclassification rate of observations from $\mathcal{C}_1$, since the misclassification rate of observations from $\mathcal{C}_2$ can be easily obtained by interchanging $n_1$ with $n_2$ and $\boldsymbol{\mu}_1$ with $\boldsymbol{\mu}_2$. The high dimensionality is modeled through its dependence on $n$, namely $p_n \to \infty$. However, we will suppress its dependence on $n$ whenever there is no confusion.

Let $\mathbf{R} = \mathbf{D}^{-1/2}\boldsymbol{\Sigma}\mathbf{D}^{-1/2}$ be the correlation matrix, and $\lambda_{\max}(\mathbf{R})$ be its largest eigenvalue, and $\boldsymbol{\alpha} \equiv (\alpha_1, \ldots, \alpha_p)' = \boldsymbol{\mu}_1 - \boldsymbol{\mu}_2$. Consider the parameter space

$$\Gamma = \left\{ (\boldsymbol{\alpha}, \boldsymbol{\Sigma}) : \boldsymbol{\alpha}'\mathbf{D}^{-1}\boldsymbol{\alpha} \geq C_p, \lambda_{\max}(\mathbf{R}) \leq b_0, \min_{1 \leq j \leq p, k=1,2} \sigma_{kj}^2 > 0 \right\},$$

where $C_p$ is a deterministic positive sequence that depends only on the dimensionality $p$, and $b_0$ is a positive constant. Note that $\boldsymbol{\alpha}'\mathbf{D}^{-1}\boldsymbol{\alpha}$ corresponds to the overall strength of signals, and the first condition $\boldsymbol{\alpha}'\mathbf{D}^{-1}\boldsymbol{\alpha} \geq C_p$ imposes a lower bound on the strength of signals. The second condition $\lambda_{\max}(\mathbf{R}) \leq b_0$ requires that the maximum eigenvalue of $\mathbf{R}$ should not exceed a positive constant. But since there are no restrictions on the smallest eigenvalue of $\mathbf{R}$, the condition number can still diverge. The third condition $\min_{1 \leq j \leq p, k=1,2} \sigma_{kj}^2 > 0$ ensures that there are no deterministic features that make classification trivial and the diagonal matrix $\mathbf{D}$ is always invertible. We will consider the asymptotic behavior of $W(\hat{\delta}, \boldsymbol{\theta})$ and $W(\hat{\delta})$.

THEOREM 1.   *Suppose that* $\log p = o(n)$, $n = o(p)$ *and* $nC_p \to \infty$. *Then:*

(i) *The classification error* $W(\delta, \boldsymbol{\theta})$ *with* $\boldsymbol{\theta} \in \Gamma$ *is bounded from above as*

$$W(\hat{\delta}, \boldsymbol{\theta}) \leq 1 - \Phi\left( \frac{[n_1 n_2/(pn)]^{1/2}\boldsymbol{\alpha}'\mathbf{D}^{-1}\boldsymbol{\alpha}(1 + o_P(1)) + \sqrt{p/(nn_1 n_2)}(n_1 - n_2)}{2\sqrt{\lambda_{\max}(\mathbf{R})}\{1 + n_1 n_2/(pn)\boldsymbol{\alpha}'\mathbf{D}^{-1}\boldsymbol{\alpha}(1 + o_P(1))\}^{1/2}} \right).$$

(ii) *Suppose* $p/(nC_p) \to 0$. *For the worst case classification error* $W(\delta)$, *we have*

$$W(\hat{\delta}) = 1 - \Phi(\tfrac{1}{2}[n_1 n_2/(pnb_0)]^{1/2}C_p\{1 + o_P(1)\}).$$

*Specifically, when* $\{\frac{n_1 n_2}{pn}\}^{1/2}C_p \to C_0$ *with* $C_0$ *a nonnegative constant, then*

$$W(\hat{\delta}) \xrightarrow{P} 1 - \Phi(C_0/(2\sqrt{b_0})).$$

*In particular, if* $C_0 = 0$, *then* $W(\hat{\delta}) \xrightarrow{P} \frac{1}{2}$.

Theorem 1 reveals the trade-off between the signal strength $C_p$ and the dimensionality, reflected in the term $C_p/\sqrt{p}$ when all features are used for classification. It states that the independence rule $\hat{\delta}$ would be no better than the random guessing due to noise accumulation, unless the signal levels are extremely high, say, $\{\frac{n}{p}\}^{1/2}C_p \geq B$ for some $B > 0$. Indeed, discrimination based on linear projections to almost all directions performs nearly the same as random guessing, as shown in the theorem below. The poor performance is caused by noise accumulation in the estimation of $\boldsymbol{\mu}_1$ and $\boldsymbol{\mu}_2$.

THEOREM 2. *Suppose that* $\mathbf{a}$ *is a p-dimensional uniformly distributed unit random vector on a $(p-1)$-dimensional sphere. Let $\lambda_1, \ldots, \lambda_p$ be the eigenvalues of the covariance matrix $\boldsymbol{\Sigma}$. Suppose $\lim_p \frac{1}{p^2} \sum_{j=1}^{p} \lambda_j^2 < \infty$ and $\lim_p \frac{1}{p} \sum_{j=1}^{p} \lambda_j = \tau$ with $\tau$ a positive constant. Moreover, assume that $p^{-1} \boldsymbol{\alpha}' \boldsymbol{\alpha} \to 0$. Then if we project all the observations onto the vector $\mathbf{a}$ and use the classifier*

$$(2.3) \qquad \hat{\delta}_{\mathbf{a}}(\mathbf{x}) = (\mathbf{a}'\mathbf{x} - \mathbf{a}'\widehat{\boldsymbol{\mu}})(\mathbf{a}'\widehat{\boldsymbol{\mu}}_1 - \mathbf{a}'\widehat{\boldsymbol{\mu}}_2),$$

*the misclassification rate of $\hat{\delta}_{\mathbf{a}}$ satisfies*

$$P(\hat{\delta}_{\mathbf{a}}(\mathbf{X}) \leq 0 | \mathbf{Y}_{ki}, \; i = 1, \ldots, n_k, \; k = 1, 2) \xrightarrow{P} \tfrac{1}{2},$$

*where the probability is taken with respect to $\mathbf{a}$ and $\mathbf{X} \in \mathcal{C}_1$.*

**3. Feature selection by two-sample $t$-test.** To extract salient features, we appeal to the two-sample $t$-test statistics. Other componentwise tests such as the rank sum test can also be used, but we do not pursue those in detail. The two-sample $t$-statistic for feature $j$ is defined as

$$(3.1) \qquad T_j = \frac{\bar{Y}_{1j} - \bar{Y}_{2j}}{\sqrt{S_{1j}^2/n_1 + S_{2j}^2/n_2}}, \qquad j = 1, \ldots, p,$$

where $\bar{Y}_{kj}$ and $S_{kj}^2$ are the same as those defined in Section 1. We work under more relaxed technical conditions: the normality assumption is not needed. Instead, we assume merely that the noise vectors $\boldsymbol{\epsilon}_{ki}, \; i = 1, \ldots, n_k$, are i.i.d. within class $\mathcal{C}_k$ with mean $\mathbf{0}$ and covariance matrix $\boldsymbol{\Sigma}_k$, and are independent between classes. The covariance matrix $\boldsymbol{\Sigma}_1$ can also differ from $\boldsymbol{\Sigma}_2$.

To show that the $t$-statistic can select all the important features with probability 1, we need the following condition.

CONDITION 1.

(a) Assume that the vector $\boldsymbol{\alpha} = \boldsymbol{\mu}_1 - \boldsymbol{\mu}_2$ is sparse and without loss of generality, only the first $s$ entries are nonzero.

(b) Suppose that $\epsilon_{kij}$ and $\epsilon_{kij}^2 - 1$ satisfy the Cramér's condition, that is, there exist constants $\nu_1$, $\nu_2$, $M_1$ and $M_2$, such that $E|\epsilon_{kij}|^m \leq m! M_1^{m-2} \nu_1/2$ and $E|\epsilon_{kij}^2 - \sigma_{kj}^2|^m \leq m! M_2^{m-2} \nu_2/2$ for all $m = 1, 2, \ldots$.

(c) Assume that the diagonal elements of both $\boldsymbol{\Sigma}_1$ and $\boldsymbol{\Sigma}_2$ are bounded away from 0.

The following theorem describes the situation under which the two-sample $t$-test can pick up all important features by choosing an appropriate critical value. Recall that $c_1 \leq n_1/n_2 \leq c_2$ and $n = n_1 + n_2$.

THEOREM 3. *Let $s$ be a sequence such that $\log(p - s) = o(n^\gamma)$ and $\log s = o(n^{1/2-\gamma}\beta_n)$ for some $\beta_n \to \infty$ and $0 < \gamma < \frac{1}{3}$. Suppose that $\min_{1 \leq j \leq s} \frac{|\alpha_j|}{\sqrt{\sigma_{1j}^2 + \sigma_{2j}^2}} = n^{-\gamma}\beta_n$. Then under Condition 1, for $x \sim cn^{\gamma/2}$ with $c$ some positive constant, we have*

$$P\left(\min_{j \leq s} |T_j| \geq x \text{ and } \max_{j > s} |T_j| < x\right) \to 1.$$

In the proof of Theorem 3, we used the moderate deviation results of the two-sample $t$-statistic [see Cao (2007) or Shao (2005)]. Theorem 3 allows the lowest signal level to decay with sample size $n$. As long as the rate of decay is not too fast and the sample size is not too small, the two-sample $t$-test can pick up all the important features with probability tending to 1.

**4. Features annealed independence rules.** We apply the independence classifier to the selected features, resulting in a Features Annealed Independence Rule (FAIR). In many applications such as tumor classification using gene expression data, we would expect that elements in the population mean difference vector $\boldsymbol{\alpha}$ are sparse: most entries are small. Thus, even if we could use $t$-test to correctly extract out all these features, the resulting choice is not necessarily optimal, since the noise accumulation can even exceed the signal accumulation for faint features. This can be seen from Theorem 1. Therefore, it is necessary to further single out the most important features that help reduce misclassification rate.

To help us select the number of features, or the critical value of the test statistic, we first consider the ideal situation that the important features are located at the first $m$ coordinates and our task is to merely select $m$ to minimize the misclassification rate. This is the case when we have the ideal information about the relative importance of features, as measured by $|\alpha_j|/\sigma_j$, say. When such an oracle information is unavailable, we will learn it from the data. In the situation that we have vague knowledge about the importance of features such as tumor classification using gene expression data, we can give high ranks to features with large $|\alpha_j|/\sigma_j$.

In the presentation below, unless otherwise specified, we assume that the two classes $\mathcal{C}_1$ and $\mathcal{C}_2$ are both from Gaussian distributions and the common covariance matrix is the identity, that is, $\boldsymbol{\Sigma}_1 = \boldsymbol{\Sigma}_2 = \mathbf{I}$. If this common covariance matrix is known, the independence classifier $\hat{\delta}$ becomes the nearest centroids classifier

$$\hat{\delta}_{\mathrm{NC}}(\mathbf{x}) = (\mathbf{x} - \hat{\boldsymbol{\mu}})'(\hat{\boldsymbol{\mu}}_1 - \hat{\boldsymbol{\mu}}_2).$$

If only the first $m$ dimensions are used in the classification, the corresponding features annealed independence classifier becomes

$$\hat{\delta}_{\mathrm{NC}}^m(\mathbf{x}) = (\mathbf{x}^m - \hat{\boldsymbol{\mu}}^m)'(\hat{\boldsymbol{\mu}}_1^m - \hat{\boldsymbol{\mu}}_2^m),$$

where the superscript $m$ means that the vector is truncated after the first $m$ entries. This is indeed the same as the nearest shrunken centroids method of Tibshirani et al. (2002).

THEOREM 4. *Consider the truncated classifier $\hat{\delta}_{NC}^{m_n}$ for a given sequence $m_n$. Suppose that $\frac{n}{\sqrt{m_n}}\sum_{j=1}^{m_n}\alpha_j^2 \to \infty$ as $m_n \to \infty$. Then the classification error of $\hat{\delta}_{NC}^{m_n}$ is*

$$W(\hat{\delta}_{NC}^{m_n}, \boldsymbol{\theta}) = 1 - \Phi\left(\frac{(1+o_P(1))\sum_{j=1}^{m_n}\alpha_j^2 + m_n(n_1-n_2)/(n_1 n_2)}{2\{(1+o_P(1))\sum_{j=1}^{m_n}\alpha_j^2 + n m_n/(n_1 n_2)\}^{1/2}}\right),$$

*where $n = n_1 + n_2$ as defined in Section 2.*

In the following, we suppress the dependence of $m$ on $n$ when there is no confusion. The above theorem reveals that the ideal choice on the number of features is

$$m_0 = \arg\max_{1 \le m \le p} \frac{[\sum_{j=1}^{m}\alpha_j^2 + m(n_1-n_2)/(n_1 n_2)]^2}{nm/(n_1 n_2) + \sum_{j=1}^{m}\alpha_j^2}.$$

It can be estimated as

$$\hat{m}_0 = \arg\max_{1 \le m \le p} \frac{[\sum_{j=1}^{m}\hat{\alpha}_j^2 + m(n_1-n_2)/(n_1 n_2)]^2}{nm/(n_1 n_2) + \sum_{j=1}^{m}\hat{\alpha}_j^2},$$

where $\hat{\alpha}_j = \hat{\mu}_{1j} - \hat{\mu}_{2j}$. The expression for $m_0$ quantifies how the signal and the noise affect the misclassification rates as the dimensionality $m$ increases. In particular, when $n_1 = n_2$, the express reduces to $m_0 = \arg\max_{1 \le m \le p} \frac{[m^{-1/2}\sum_{j=1}^{m}\alpha_j^2]^2}{2/n + \sum_{j=1}^{m}\alpha_j^2/m}$. The term $m^{-1/2}\sum_{j=1}^{m}\alpha_j^2$ reflects the trade-off between the signal and noise as dimensionality $m$ increases.

The good performance of the classifier $\hat{\delta}_{NC}^m$ depends on the assumption that the largest entries of $\boldsymbol{\alpha}$ cluster at the first $m$ dimensions. An ideal version of the classifier $\hat{\delta}_{NC}$ is to select a subset $\mathcal{A} = \{j : |\alpha_j| > a\}$ and use this subset to construct independence classifier. Let $m$ be the number of elements in $\mathcal{A}$. The oracle classifier can be written as

$$\hat{\delta}_{orc}(\mathbf{x}) = \sum_{j=1}^{p} \hat{\alpha}_j(x_j - \hat{\mu}_j) 1_{\{|\alpha_j| > a\}}.$$

The misclassification rate is approximately

(4.1) $$1 - \Phi\left(\frac{\sum_{j \in \mathcal{A}}\alpha_j^2 + m(n_1-n_2)/(n_1 n_2)}{2\{nm/(n_1 n_2) + \sum_{j \in \mathcal{A}}\alpha_j^2\}^{1/2}}\right),$$

when $\frac{n}{\sqrt{m}}\sum_{j \in \mathcal{A}}\alpha_j^2 \to \infty$ and $m \to \infty$. This is straightforward from Theorem 4. In practice, we do not have such an oracle, and selecting the subset $\mathcal{A}$ is difficult. A simple procedure is to use the feature annealed independence rule based on the hard thresholding:

$$\hat{\delta}^b_{\text{FAIR}}(\mathbf{x}) = \sum_{j=1}^p \hat{\alpha}_j(x_j - \hat{\mu}_j)1_{\{|\hat{\alpha}_j| > b\}}.$$

We study the classification error of FAIR and the impact of the threshold $b$ on the classification result in the following theorem.

THEOREM 5. *Suppose that* $\max_{j \in \mathcal{A}^c} |\alpha_j| < b_n$ *and* $\log(p - m)/[n(b_n - \max_{j \in \mathcal{A}^c} |\alpha_j|)^2] \to 0$ *with* $m = |\mathcal{A}|$. *Moreover, assume that* $\frac{n}{\sqrt{m}} \sum_{j \in \mathcal{A}} \alpha_j^2 \to \infty$ *and* $\sum_{j \in \mathcal{A}} |\alpha_j|/[\sqrt{n} \sum_{j \in \mathcal{A}} \alpha_j^2] \to 0$. *Then*

$$W(\hat{\delta}^{b_n}_{\text{FAIR}}, \boldsymbol{\theta}) \leq 1 - \Phi\bigg(\frac{(1 + o_P(1)) \sum_{j \in \mathcal{A}} \alpha_j^2 + nm(n_1 n_2)^{-1} - mb_n^2}{2\{(1 + o_P(1)) \sum_{j \in \mathcal{A}} \alpha_j^2 + nm(n_1 n_2)^{-1}\}^{1/2}}\bigg).$$

Notice that the upper bound of $W(\hat{\delta}^{b_n}_{\text{FAIR}}, \boldsymbol{\theta})$ in Theorem 5 is greater than the classification error in Theorem 4, and the magnitude of difference depends on $mb_n^2$. This is expected as estimating the set $\mathcal{A}$ increases the classification error. These results are similar to those in Fan (1996) for high-dimensional hypothesis testing.

When the common covariance matrix is different from the identity, FAIR takes a slightly different form to adapt to the unknown componentwise variance:

$$(4.2) \qquad \hat{\delta}_{\text{FAIR}}(\mathbf{x}) = \sum_{j=1}^p \hat{\alpha}_j(x_j - \hat{\mu}_j)/\hat{\sigma}_j^2 1_{\{\sqrt{n/(n_1 n_2)}|T_j| > b\}},$$

where $T_j$ is the two-sample $t$-statistic. It is clear from (4.2) that FAIR works the same way as that we first sort the features by the absolute values of their $t$-statistics in the descending order, and then take out the first $m$ features to classify the data. The number of features can be selected by minimizing the upper bound of the classification error given in Theorem 1. The optimal $m$ in this sense is

$$m_1 = \arg\max_{1 \leq m \leq p} \frac{1}{\lambda^m_{\max}} \frac{[\sum_{j=1}^m \alpha_j^2/\sigma_j^2 + m(1/n_2 - 1/n_1)]^2}{nm/(n_1 n_2) + \sum_{j=1}^m \alpha_j^2/\sigma_j^2},$$

where $\lambda^m_{\max}$ is the largest eigenvalue of the correlation matrix $\mathbf{R}^m$ of the truncated observations. It can be estimated from the samples:

$$(4.3) \qquad \begin{aligned} \hat{m}_1 &= \arg\max_{1 \leq m \leq p} \frac{1}{\hat{\lambda}^m_{\max}} \frac{[\sum_{j=1}^m \hat{\alpha}_j^2/\hat{\sigma}_j^2 + m(1/n_2 - 1/n_1)]^2}{nm/(n_1 n_2) + \sum_{j=1}^m \hat{\alpha}_j^2/\hat{\sigma}_j^2} \\ &= \arg\max_{1 \leq m \leq p} \frac{1}{\hat{\lambda}^m_{\max}} \frac{n[\sum_{j=1}^m T_j^2 + m(n_1 - n_2)/n]^2}{mn_1 n_2 + n_1 n_2 \sum_{j=1}^m T_j^2}. \end{aligned}$$

Note that the factor $\lambda_{\max}^m$ in (4.3) increases with $m$, which makes $\hat{m}_1$ usually smaller than $\hat{m}_0$.

**5. Numerical studies.** In this section we use a simulation study and three real data analyses to illustrate our theoretical results and to verify the performance of our newly proposed classifier FAIR.

5.1. *Simulation study.* We first introduce the model. The covariance matrices $\boldsymbol{\Sigma}_1$ and $\boldsymbol{\Sigma}_2$ for the two classes are chosen to be the same. For the distribution of the error $\epsilon_{ij}$ in (2.1), we use the same model as that in Fan, Hall and Yao (2006). Specifically, features are divided into three groups. Within each group, features share one unobservable common factor with different factor loadings. In addition, there is an unobservable common factor among all the features across three groups. For simplicity, we assume that the number of features $p$ is a multiple of 3. Let $Z_{ij}$ be a sequence of independent $N(0,1)$ random variables, and $\chi_{ij}^2$ be a sequence of independent random variables of the same distribution as $(\chi_d^2 - d)/\sqrt{2d}$ with $\chi_d^2$ the Chi-square distribution with degrees of freedom $d$. In the simulation we set $d = 6$.

Let $\{a_j\}$ and $\{b_j\}$ be factor loading coefficients. Then the error in (2.1) is defined as

$$\epsilon_{ij} = \frac{Z_{ij} + a_{1j}\chi_{1i} + a_{2j}\chi_{2i} + a_{3j}\chi_{3i} + b_j\chi_{4i}}{(1 + a_{1j}^2 + a_{2j}^2 + a_{3j}^2 + b_j^2)^{1/2}}, \qquad i = 1, \ldots, n_k, \ j = 1, \ldots, p,$$

where $a_{ij} = 0$ except that $a_{1j} = a_j$ for $j = 1, \ldots, p/3$, $a_{2j} = a_j$ for $j = (p/3) + 1, \ldots, 2p/3$, and $a_{3j} = a_j$ for $j = (2p/3) + 1, \ldots, p$. Therefore, $E\epsilon_{ij} = 0$ and $\operatorname{var}(\epsilon_{ij}) = 1$, and in general, within group correlation is greater than the between group correlation. The factor loadings $a_j$ and $b_j$ are independently generated from uniform distributions $U(0, 0.4)$ and $U(0, 0.2)$. The mean vector $\boldsymbol{\mu}_1$ for class $\mathcal{C}_1$ is taken from a realization of the mixture of a point mass at 0 and a double-exponential distribution:

$$(1 - c)\delta_0 + \tfrac{1}{2}c\exp(-2|x|),$$

where $c \in (0,1)$ is a constant. In the simulation, we set $p = 4500$ and $c = 0.02$. In other words, there are around 90 signal features on an average, many of which are weak signals. Without loss of generality, $\boldsymbol{\mu}_2$ is set to be 0. Figure 1 shows the true mean difference vector $\boldsymbol{\alpha}$, which is fixed across all simulations. It is clear that there are only very few features with signal levels exceeding 1 standard deviation of the noise.

With the parameters and model above, for each simulation, we generate $n_1 = 30$ training data from class $\mathcal{C}_1$ and $n_2 = 30$ training data from class $\mathcal{C}_2$. In addition, separate 200 samples are generated from each of the two classes in each simulation, and these 400 vectors are used as test samples. We apply our newly proposed classifier FAIR to the simulated data. Specifically, for each feature, the $t$-test statistic in (3.1) is calculated using the training sample. Then the features are sorted in the decreasing order of the absolute values of their $t$-statistics. We then examine the impact of the number of features $m$ on the misclassification rate. In each simulation, with $m$ ranging from 1 to 4500, we construct the feature annealed independence classifiers using the training samples, and then apply these classifiers to the 400 test samples. The classification errors are compared to those of the independence rule with the oracle ordering information, which is constructed by repeating the above procedure except that in the first step the features are ordered by their true signal levels, $|\boldsymbol{\alpha}|$, instead of by their $t$-statistics.

The above procedure is repeated 100 times, and averages and standard errors of the misclassification rates (based on 400 test samples in each simulation) are calculated across the 100 simulations. Note that the average of the 100 misclassification rates is indeed computed based on $100 \times 400$ testing samples.

Figure 2 depicts the misclassification rate as a function of the number of features $m$. The solid curves represent the average of classification rates across the 100 simulations, and the corresponding dashed curves are 2 standard errors (i.e., the standard deviation of 100 misclassification rates divided by 10) away from the solid one. The misclassification rates using the first 80 features in Figure 2(a) are zoomed in Figure 2(b). Figures 2(c) and 2(d) are the same as 2(a) and 2(b) except that the features are arranged in the decreasing order of $|\boldsymbol{\alpha}|$, that is, the results are based on the oracle-assisted feature annealed independence classifier. From these plots we see that the

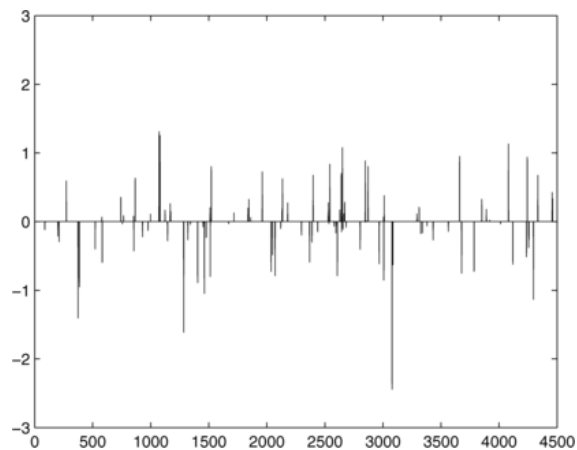

FIG. 1. *True mean difference vector $\boldsymbol{\alpha}$. x-axis represents the dimensionality, and y-axis shows the values of corresponding entries of $\boldsymbol{\alpha}$.*

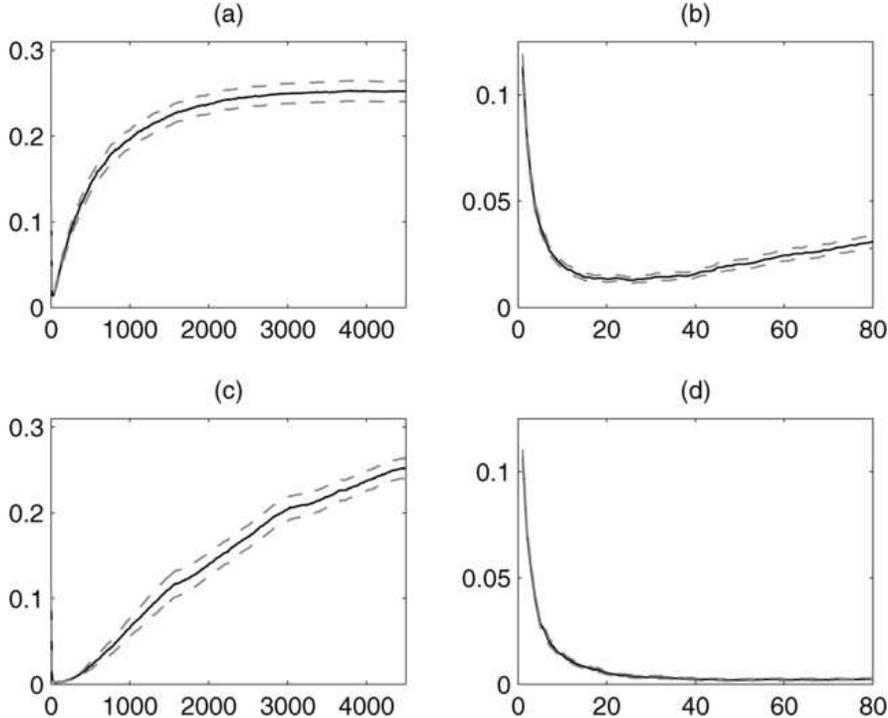

Fig. 2.    *Number of features versus misclassification rates. The solid curves represent the averages of classification errors across 100 simulations. The dashed curves are 2 standard errors away from the solid curves. The x-axis represents the number of features used in the classification, and the y-axis shows the misclassification rates. (a) The features are ordered in a way such that the corresponding t-statistics are decreasing in absolute values. (b) The amplified plot of the first 80 values of x-axis in plot (*a*). (c) The same as (*a*) except that the features are arranged in a way such that the corresponding true mean differences are decreasing in absolute values. (d) The amplified plot of the first 80 values of x-axis in plot (*c*).*

classification results of FAIR are close to those of the oracle-assisted independence classifier. Moreover, as the dimensionality $m$ grows, the misclassification rate increases steadily due to the noise accumulation. When all the features are included, that is, $m = 4500$, the misclassification rate is 0.2522, whereas the minimum classification errors are 0.0128 in plot 2(b) and 0.0020 in plot 2(d). These results are consistent with Theorem 1. We also tried to decrease the signal levels, that is, the mean of the double exponential distribution, or to increase the dimensionality $p$, and found that the classification error tend to 0.5 when all the dimensions are included. Comparing Figures 2(a) and 2(b) to Figures 2(c) and 2(d), we see that the features ordered by $t$-statistics has higher misclassification rates than those ordered by the oracle. Also, using $t$-statistics results in larger minimum classification errors [see plots 2(b) and 2(d)], but the differences are not very large.

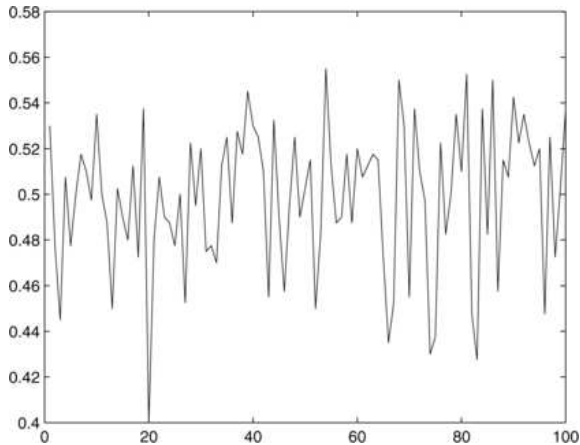

Fig. 3. *Classification errors of the independence rule based on projected samples onto randomly chosen directions over 100 simulations.*

Figure 3 shows the classification errors of the independence rule based on projected samples onto randomly chosen directions across 100 simulations. Specifically, in each of the simulations in Figure 2, we generate a direction vector **a** randomly from the $(p-1)$-dimensional unit sphere, then project all the data in that simulation onto the direction **a**, and finally apply the Fisher discriminant to the projected data [see (2.3)]. The average of these misclassification rates is 0.4986 and the corresponding standard deviation is 0.0318. These results are consistent with our Theorem 2.

Finally, we examine the effectiveness of our proposed method (4.3) for selecting features in FAIR. In each of the 100 simulations, we apply (4.3) to choose the number of features and compute the resulting misclassification rate based on 400 test samples. We also use the nearest shrunken centroids of Tibshirani et al. (2002) to select the important features. Figure 4 summarizes these results. The thin curves correspond to the nearest shrunken centroids method, and the thick curves correspond to FAIR. Figure 4(a) presents the number of features calculated from these two methods, and Figure 4(b) shows the corresponding misclassification rates. For our newly proposed classifier FAIR, the average of the optimal number of features over 100 simulations is 29.71, which is very close to the smallest number of features with the minimum misclassification rate in Figure 2(d). The misclassification rates of FAIR in Figure 4(b) have average 0.0154 and standard deviation 0.0085, indicating the outstanding performance of FAIR. Nearest shrunken centroids method is unstable in selecting features. Over the 100 simulations, there are several realizations in which it chooses plenty of fea-tures. We truncated Figure 4 to make it easier to view. The average number of features chosen by the nearest shrunken centroids is 28.43, and the average classification error is 0.0216, with corresponding standard deviation 0.0179. It is clear that nearest shrunken centroids method tends to choose less features than FAIR, but the misclassification rates are larger.

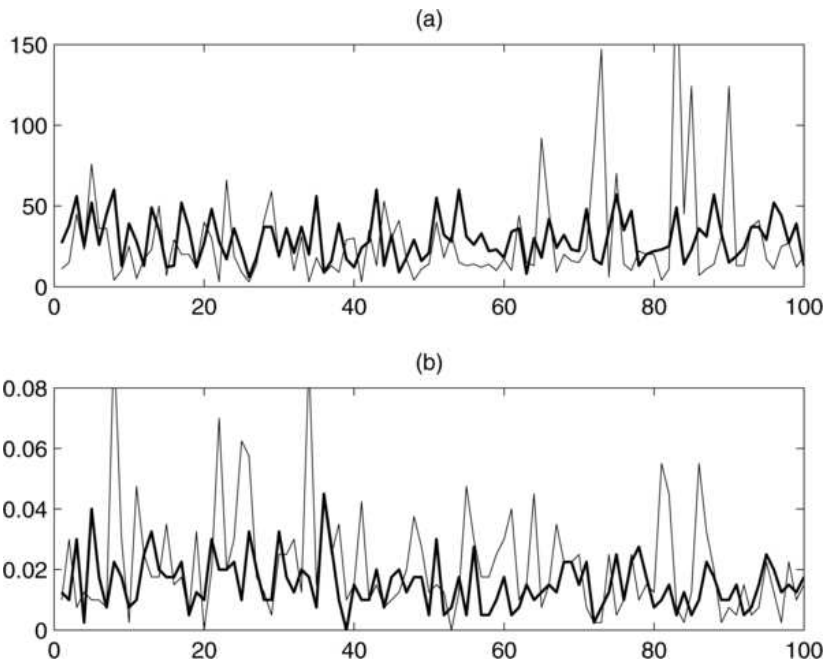

FIG. 4. *The thick curves correspond to FAIR, while the thin curves correspond to the nearest shrunken centroids method.* (a) *The numbers of features chosen by* (4.3) *and by the nearest shrunken centroids method over 100 simulations.* (b) *Corresponding classification errors based on the optimal number of features chosen in* (a) *over 100 simulations.*

### 5.2. *Real data analysis.*

5.2.1. *Leukemia data.* Leukemia data from high-density Affymetrix oligonucleotide arrays were previously analyzed in Golub et al. (1999), and are available at http://www.broad.mit.edu/cgi-bin/cancer/datasets.cgi. There are 7129 genes and 72 samples coming from two classes: 47 in class ALL (acute lymphocytic leukemia) and 25 in class AML (acute mylogenous leukemia). Among these 72 samples, 38 (27 in class ALL and 11 in class AML) are set to be training samples and 34 (20 in class ALL and 14 in class AML) are set as test samples.

Before classification, we standardize each sample to zero mean and unit variance as done by Dudoit, Fridlyand and Speed (2002). The classification results from the nearest shrunken centroids (NSC hereafter) method and FAIR are shown in Table 1. The nearest shrunken centroids method picks up 21 genes and makes 1 training error and 3 test errors, while our method chooses 11 genes and makes 1 training error and 1 test error. Tibshirani et al. (2002) proposed and applied the nearest shrunken centroids method to the unstandardized Leukemia dataset. They chose 21 genes and made 1 training error and 2 test errors. Our results are still superior to theirs.

To further evaluate the performance of the two classifiers, we randomly split the 72 samples into training and test sets. Specifically, we set approximately $100\gamma\%$ of the observations from class ALL and $100\gamma\%$ of the observations from class AML as training samples, and the rest as test samples. FAIR and NSC are applied to the training data, and their performances are evaluated by the test samples. The above procedure is repeated 100 times for $\gamma = 0.4$, 0.5 and 0.6, respectively, and the distributions of test errors of FAIR, NSC and the independence rule without feature selection are summarized in Figure 5. In each of the splits, we also calculated the difference of test errors between NSC and FAIR, that is, the test error of FAIR minus that of NSC, and the distribution is summarized in Figure 5. The top panel of Figure 6 shows the number of features selected by FAIR and NSC for $\gamma = 0.4$. The results for the other two values of $\gamma$ are similar so we do not present here to save the space. From these figures we can see that the performance of independence rule improves significantly after feature selection. The classification errors of NSC and FAIR are approximately the same. As we have already noticed in the simulation study, NSC is not good with feature selection, that is, the number of features selected by NSC is very large and unstable, while the number of features selected by FAIR is quite reasonable and stable over different random splits. Clearly, the independent rule without feature selection performs poorly.

5.2.2. *Lung cancer data.* We evaluate our method by classifying between malignant pleural mesothelioma (MPM) and adenocarcinoma (ADCA) of the lung. Lung cancer data were analyzed by Gordon et al. (2002) and are available at **http://www.chestsurg.org**. There are 181 tissue samples (31 MPM and 150 ADCA). The training set contains 32 of them, with 16 from MPM and 16 from ADCA. The rest 149 samples are used for testing (15 from MPM and 134 from ADCA). Each sample is described by 12533 genes.

TABLE 1
*Classification errors of Leukemia dataset*

| Method | Training error | Test error | No. of selected genes |
|---|---|---|---|
| Nearest shrunken centroids | 1/38 | 3/34 | 21 |
| FAIR | 1/38 | 1/34 | 11 |

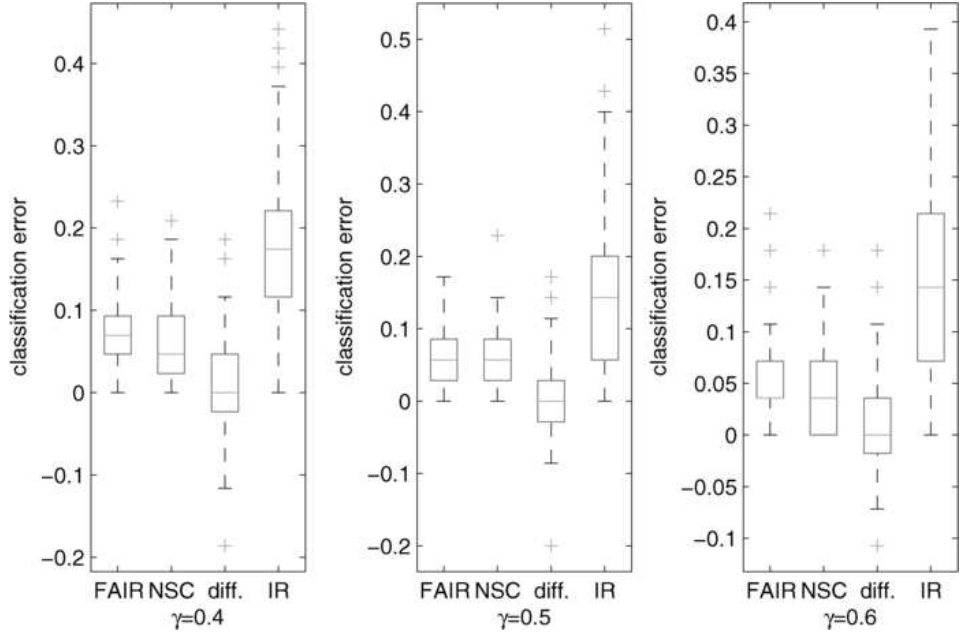

Fig. 5.  *Leukemia data. Boxplots of test errors of FAIR, NSC and the independence rule without feature selection over 100 random splits of 72 samples, where $100\gamma\%$ of the samples from both classes are set as training samples. The three plots from left to right correspond to $\gamma = 0.4, 0.5$ and 0.6, respectively. In each boxplot above, "FAIR" refers to the test errors of the feature annealed independent rule; "NSC" corresponds to the test errors of nearest shrunken centroids method; "diff." means the difference of the test errors of FAIR and those of NSC; and "IR" corresponds the test errors of independence rule without feature selection.*

As in the Leukemia dataset, we first standardize the data to zero mean and unit variance, and then apply the two classification methods to the standardized dataset. Classification results are summarized in Table 2. Although FAIR uses 5 more genes than the nearest shrunken centroids method, it has better classification results: both methods perfectly classify the training samples, while our classification procedure has smaller test error.

We follow the same procedure as that in Leukemia example to randomly split the 181 samples into training and test sets. FAIR and NSC are applied to the training data, and the test errors are calculated using the test data. The procedure is repeated 100 times with $\gamma = 0.4, 0.5$ and 0.6, respectively, and the test error distributions of FAIR, NSC and the independence rule without feature selection can be found in Figure 7. We also present the difference of the test errors between FAIR and NSC in Figure 7. The numbers of features used by FAIR and NSC with $\gamma = 0.5$ are shown in the middle panel of Figure 6. Figure 7 shows again that feature selection is very important in high-dimensional classification. The performance of FAIR is close to NSC in terms of classification error (Figure 7), but FAIR is stable in feature selection, as shown in the middle panel of Figure 6. One possible reason of Figure 7 might be that the signal strength in this Lung cancer dataset is relatively weak, and more features are needed to obtain the optimal performance. However, the estimate of the largest eigenvalue is not accurate anymore when the number of features is large, which results in inaccurate estimates of $m_1$ in (4.3).

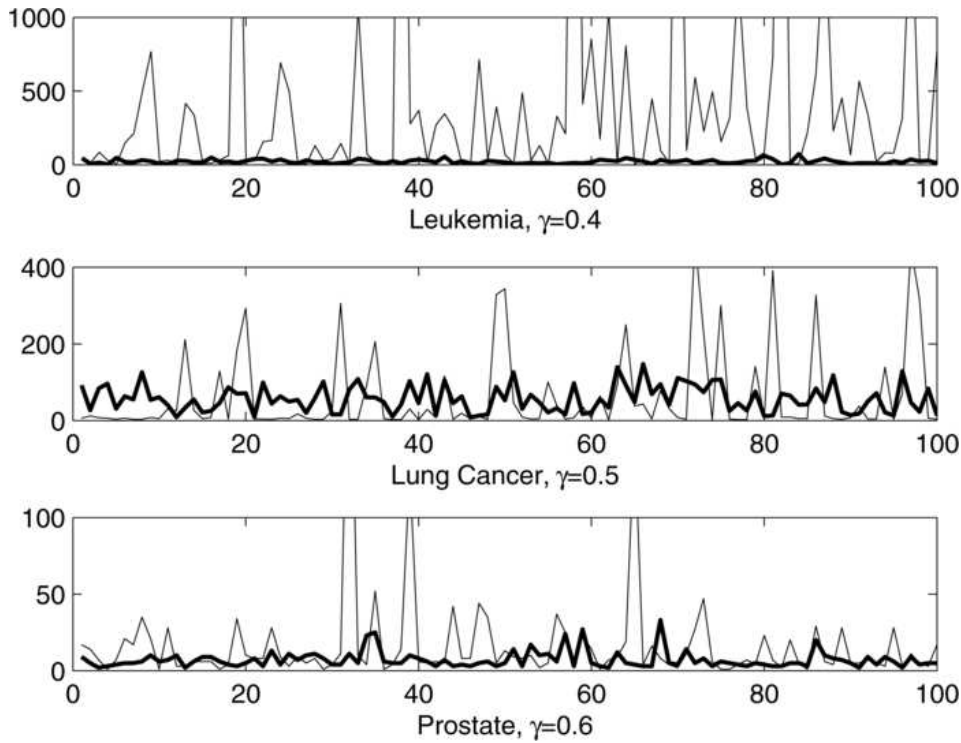

Fig. 6. *Leukemia, Lung cancer and Prostate datasets. The number of features selected by FAIR and NSC over 100 random splits of the total samples. In each split, $100\gamma\%$ of the samples from both class are set as training samples, and the rest are used as test samples. The three plots from top to bottom correspond to the Leukemia data with $\gamma = 0.4$, the Lung cancer data with $\gamma = 0.5$ and the Prostate cancer data with $\gamma = 0.6$, respectively. The thin curves show the results from NSC, and the thick curves correspond to FAIR. The plots are truncated to make them easy to view.*

5.2.3. *Prostate cancer data.* The last example uses the prostate cancer data studied in Singh et al. (2002). The dataset is available at

Table 2
*Classification errors of Lung cancer data*

| Method | Training error | Test error | No. of selected genes |
|---|---|---|---|
| Nearest shrunken centroids | 0/32 | 11/149 | 26 |
| FAIR | 0/32 | 7/149 | 31 |

Table 3
*Classification errors of Prostate cancer dataset*

| Method | Training error | Test error | No. of selected genes |
|---|---|---|---|
| Nearest shrunken centroids | 8/102 | 9/34 | 6 |
| FAIR | 10/102 | 9/34 | 2 |

http://www.broad.mit.edu/cgi-bin/cancer/datasets.cgi. The training dataset contains 102 patient samples, 52 of which (labeled as "tumor") are prostate tumor samples and 50 of which (labeled as "Normal") are prostate samples. There are around 12600 genes. An independent set of test samples is from a different experiment and has 25 tumor and 9 normal samples.

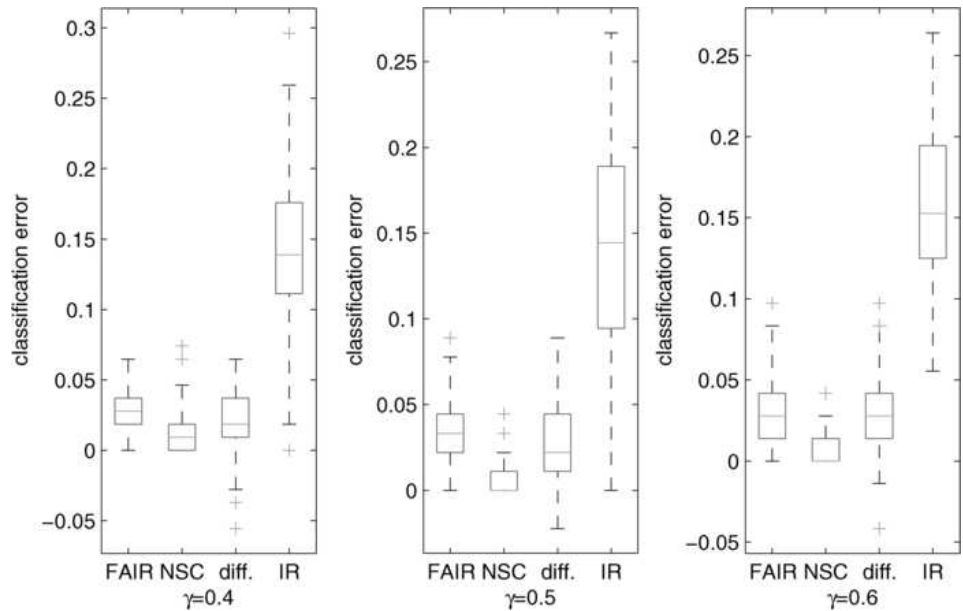

Fig. 7. *Lung cancer data. The same as Figure 5 except that the dataset is different.*

We preprocess the data by standardizing the gene expression data as before. The classification results are summarized in Table 3. We make the same test error as and a bit larger training error than the nearest shrunken centroids method, but the number of selected genes we use is much less.

The samples are randomly split into training and test sets in the same way as before, the test errors are calculated, and the number of features used by these two methods are recorded. Figure 8 shows the test errors of FAIR, NSC and the independence rule without feature selection, and the difference of the test errors of FAIR and NSC. The bottom panel of Figure 6 presents the numbers of features used by FAIR and NSC in each random split for $\gamma = 0.6$. As we mentioned before, the plots for $\gamma = 0.4$ and $0.5$ are similar so we omit them in the paper. The performance of FAIR is better than that of NSC both in terms of classification error and in terms of the selection of features. The good performance of FAIR might be caused by the strong signal level of few features in this dataset. Due to the strong signal level, FAIR can attain the optimal performance with small number of features. Thus, the estimate of $m_1$ in (4.3) is accurate and hence the actual performance of FAIR is good.

**6. Conclusion.** This paper studies the impact of high dimensionality on classifications. To illustrate the idea, we have considered the independence classification rule, which avoids the difficulty of estimating large covariance matrix and the diverging condition number frequently associated with the large covariance matrix. When only a subset of the features capture the characteristics of two groups, classification using all dimensions would intrinsically classify the noises. We prove that classification based on linear projections onto almost all directions performs nearly the same as random guessing. Hence, it is necessary to choose direction vectors which put more weight on important features.

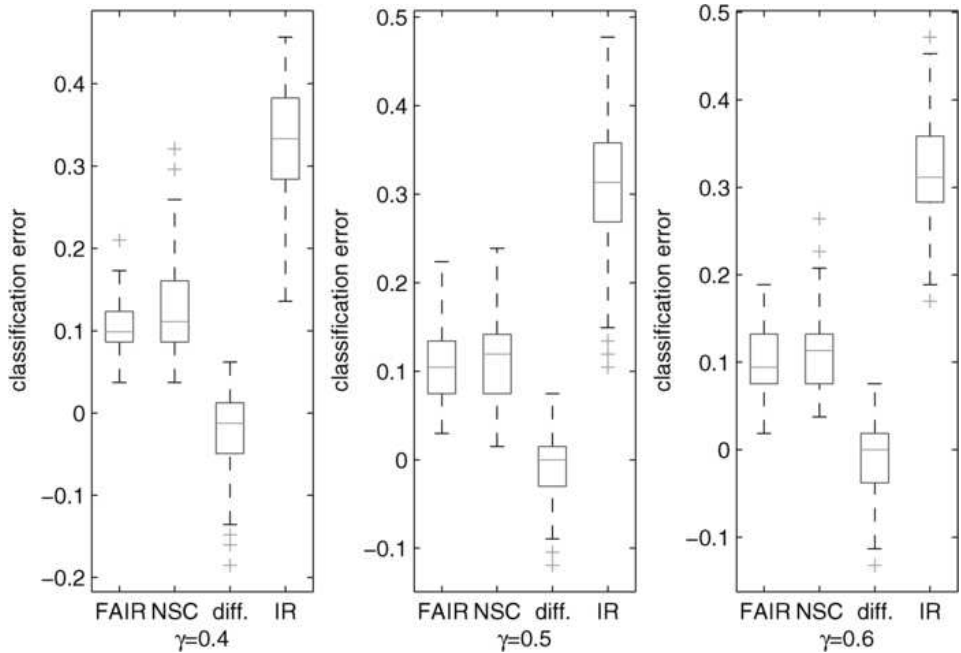

Fig. 8. *Prostate cancer data. The same as Figure 5 except that the dataset is different.*

The two-sample $t$-test can be used to choose the important features. We have shown that under mild conditions, the two-sample $t$-test can select all the important features with probability one. The features annealed independence rule using hard thresholding, FAIR, is proposed, with the number of features selected by a data-driven rule. An upper bound of the classification error of FAIR is explicitly given. We also give suggestions on the optimal number of features used in classification. Simulation studies and real data analysis support our theoretical results convincingly.

## APPENDIX

PROOF OF THEOREM 1. For $\boldsymbol{\theta} \in \Gamma$, $\Psi$ defined in (2.2) can be bounded as

$$
(A.1) \qquad \Psi \geq \frac{(\boldsymbol{\mu}_1 - \widehat{\boldsymbol{\mu}})' \widehat{\mathbf{D}}^{-1} (\widehat{\boldsymbol{\mu}}_1 - \widehat{\boldsymbol{\mu}}_2)}{\sqrt{b_0 (\widehat{\boldsymbol{\mu}}_1 - \widehat{\boldsymbol{\mu}}_2)' \widehat{\mathbf{D}}^{-1} \mathbf{D} \widehat{\mathbf{D}}^{-1} (\widehat{\boldsymbol{\mu}}_1 - \widehat{\boldsymbol{\mu}}_2)}},
$$

where we have used the assumption that $\lambda_{\max}(\mathbf{R}) \leq b_0$. Denote by

$$
\widetilde{\Psi} = \frac{(\boldsymbol{\mu}_1 - \widehat{\boldsymbol{\mu}})' \widehat{\mathbf{D}}^{-1} (\widehat{\boldsymbol{\mu}}_1 - \widehat{\boldsymbol{\mu}}_2)}{\sqrt{(\widehat{\boldsymbol{\mu}}_1 - \widehat{\boldsymbol{\mu}}_2)' \widehat{\mathbf{D}}^{-1} \mathbf{D} \widehat{\mathbf{D}}^{-1} (\widehat{\boldsymbol{\mu}}_1 - \widehat{\boldsymbol{\mu}}_2)}}.
$$

We next study the asymptotic behavior of $\widetilde{\Psi}$.

Since Condition 1(b) in Section 3 is satisfied automatically for normal distribution, by Lemma A.2 below we have $\widehat{\mathbf{D}} = \mathbf{D}(1 + o_P(1))$, where $o_P(1)$ holds uniformly across all diagonal elements. Thus, the right-hand side of (A.1) can be written as

$$
\frac{1}{\sqrt{b_0}} \frac{(\boldsymbol{\mu}_1 - \widehat{\boldsymbol{\mu}})' \widehat{\mathbf{D}}^{-1} (\widehat{\boldsymbol{\mu}}_1 - \widehat{\boldsymbol{\mu}}_2)}{\sqrt{(\widehat{\boldsymbol{\mu}}_1 - \widehat{\boldsymbol{\mu}}_2)' \mathbf{D}^{-1} (\widehat{\boldsymbol{\mu}}_1 - \widehat{\boldsymbol{\mu}}_2)}} (1 + o_P(1)).
$$

We first consider the denominator. Notice that it can be decomposed as

$$
(\widehat{\boldsymbol{\mu}}_1 - \widehat{\boldsymbol{\mu}}_2)' \mathbf{D}^{-1} (\widehat{\boldsymbol{\mu}}_1 - \widehat{\boldsymbol{\mu}}_2)
$$
$$
= \boldsymbol{\alpha}' \mathbf{D}^{-1} \boldsymbol{\alpha} + 2 \sum \alpha_j \frac{\hat{\epsilon}_{1j} - \hat{\epsilon}_{2j}}{\sigma_j^2} + \sum \frac{(\hat{\epsilon}_{1j} - \hat{\epsilon}_{2j})^2}{\sigma_j^2}
$$

(A.2)

$$= \boldsymbol{\alpha}' \mathbf{D}^{-1} \boldsymbol{\alpha} + 2 \sum \frac{\alpha_j}{\sigma_j^2} (\hat{\epsilon}_{1j} - \hat{\epsilon}_{2j}) + \sum \frac{(\hat{\epsilon}_{1j} - \hat{\epsilon}_{2j})^2}{\sigma_j^2}$$

$$\equiv \boldsymbol{\alpha}' \mathbf{D}^{-1} \boldsymbol{\alpha} + 2(1 + o_P(1)) I_1 + I_2,$$

where $\sigma_j^2$ is the $j$th diagonal entry of $\mathbf{D}$, $\hat{\sigma}_j^2$ is the $j$th diagonal entry of $\widehat{\mathbf{D}}$, and $\hat{\epsilon}_{kj} = \sum_{i=1}^{n_k} \epsilon_{kij}/n_k$, $k = 1, 2$. Notice that $\hat{\boldsymbol{\epsilon}}_1 - \hat{\boldsymbol{\epsilon}}_2 \sim N(0, \frac{n}{n_1 n_2} \boldsymbol{\Sigma})$. By singular value decomposition we have

$$\mathbf{R} = \mathbf{Q}_R \mathbf{V}_R \mathbf{Q}_R',$$

where $\mathbf{Q}_R$ is orthogonal matrix and $\mathbf{V}_R = \text{diag}\{\lambda_{R,1}, \ldots, \lambda_{R,p}\}$ be the eigenvalues of the correlation matrix $\mathbf{R}$. Define $\tilde{\boldsymbol{\epsilon}} = \sqrt{n_1 n_2/n} \mathbf{V}_R^{-1/2} \mathbf{Q}_R' \mathbf{D}^{-1/2} (\hat{\boldsymbol{\epsilon}}_1 - \hat{\boldsymbol{\epsilon}}_2)$, then $\tilde{\boldsymbol{\epsilon}} \sim N(0, \mathbf{I})$. Hence,

$$I_2 = (\hat{\boldsymbol{\epsilon}}_1 - \hat{\boldsymbol{\epsilon}}_2)' \mathbf{D}^{-1} (\hat{\boldsymbol{\epsilon}}_1 - \hat{\boldsymbol{\epsilon}}_2) = \frac{n}{n_1 n_2} \tilde{\boldsymbol{\epsilon}}' \mathbf{V}_R \tilde{\boldsymbol{\epsilon}}.$$

Since $\sum_{i=1}^p \lambda_{R,i} = p$ and $\lambda_{R,i} \geq 0$ for all $i = 1, \ldots, p$, we have $\frac{1}{p^2} \sum_{i=1}^p \lambda_{R,i}^2 < \infty$. By the weak law of large number we have

(A.3)           $$n_1 n_2 I_2 / [pn] \xrightarrow{P} 1 \qquad \text{as } n \to \infty, p \to \infty.$$

Next, we consider $I_1$. Note that $I_1$ has the distribution $I_1 \sim N(0, \frac{n}{n_1 n_2} \boldsymbol{\alpha}' \mathbf{D}^{-1} \times \boldsymbol{\Sigma} \mathbf{D}^{-1} \boldsymbol{\alpha})$. Since $\lambda_{\max} \leq b_0$, $n \boldsymbol{\alpha}' \mathbf{D}^{-1} \boldsymbol{\alpha} \geq n C_p \to \infty$ and

$$\boldsymbol{\alpha}' \mathbf{D}^{-1} \boldsymbol{\Sigma} \mathbf{D}^{-1} \boldsymbol{\alpha} = \boldsymbol{\alpha}' \mathbf{D}^{-1/2} \mathbf{R} \mathbf{D}^{-1/2} \boldsymbol{\alpha} \leq \lambda_{\max}(R) \boldsymbol{\alpha}' \mathbf{D}^{-1} \boldsymbol{\alpha},$$

we have $I_1 = \boldsymbol{\alpha}' \mathbf{D}^{-1} \boldsymbol{\alpha} o_P(1)$. This together with (A.2) and (A.3) yields

(A.4)    $$\frac{n_1 n_2}{pn} (\widehat{\boldsymbol{\mu}}_1 - \widehat{\boldsymbol{\mu}}_2)' \widehat{\mathbf{D}}^{-1} (\widehat{\boldsymbol{\mu}}_1 - \widehat{\boldsymbol{\mu}}_2) = 1 + \frac{n_1 n_2}{pn} \sum_{j=1}^p \frac{\alpha_j^2}{\sigma_j^2} (1 + o_P(1)).$$

Now, we consider the numerator. It can be decomposed as

$$(\boldsymbol{\mu}_1 - \widehat{\boldsymbol{\mu}})' \widehat{\mathbf{D}}^{-1} (\widehat{\boldsymbol{\mu}}_1 - \widehat{\boldsymbol{\mu}}_2)$$

$$= \frac{1}{2} \boldsymbol{\alpha}' \widehat{\mathbf{D}}^{-1} \boldsymbol{\alpha} - \sum \frac{\alpha_j}{\hat{\sigma}_j^2} (\hat{\epsilon}_{2j}) - \frac{1}{2}(1 + o_P(1)) \sum \hat{\epsilon}_{1j}^2 / \sigma_j^2$$

$$+ \frac{1}{2}(1 + o_P(1)) \sum \hat{\epsilon}_{2j}^2 / \sigma_j^2$$

$$\equiv \frac{1}{2} \boldsymbol{\alpha}' \mathbf{D}^{-1} \boldsymbol{\alpha}(1 + o_P(1)) - I_3 - \frac{1}{2}(1 + o_P(1)) I_4 + \frac{1}{2}(1 + o_P(1)) I_5.$$

Denote by $\tilde{I}_3 = \sum \frac{\alpha_j}{\sigma_j^2} (\hat{\epsilon}_{2j})$. Note that

$$\max \left| \frac{\alpha_j}{\sigma_j^2} (\hat{\epsilon}_{2j}) - \frac{\alpha_j}{\hat{\sigma}_j^2} (\hat{\epsilon}_{2j}) \right| \leq \max \left| \frac{\sigma_j^2}{\hat{\sigma}_j^2} - 1 \right| \max \left| \frac{\alpha_j}{\sigma_j^2} \hat{\epsilon}_{2j} \right|$$

(A.5)
$$= o_P(1) \max \left| \frac{\alpha_j}{\sigma_j^2} \hat{\epsilon}_{2j} \right|.$$

Define $F_j = \sqrt{n_2} \frac{\alpha_j}{\sigma_j^2} \hat{\epsilon}_{2j} / \boldsymbol{\alpha}' \mathbf{D}^{-1} \boldsymbol{\alpha}$, then $\sigma_{F_j}^2 \equiv \text{var}(F_j) \leq 1$ for all $j$. For the normal distribution, we have the following tail probability inequality:

$$1 - \Phi(x) \leq \frac{1}{\sqrt{2\pi}} \frac{1}{x} e^{-x^2/2}.$$

Since $F_j \sim N(0, \sigma_{F_j}^2)$, by the above inequality we have

$$P(|F_j| \geq x) \leq 2 \exp\left\{ -\frac{x^2}{2C} \right\}$$

with $C$ some positive constant, for all $x > 0$ and $j = 1, \ldots, p$. By Lemma 2.2.10 of van der Vaart and Wellner [(1996), page 102], we have

$$(\boldsymbol{\alpha}' \mathbf{D}^{-1} \boldsymbol{\alpha})^{-1} E \max \left| \frac{\alpha_j}{\sigma_j^2} \hat{\epsilon}_{2j} \right| = n_2^{-1} E \max_{j \leq p} |F_j| \leq K \sqrt{C \log(p+1)/n_2} \xrightarrow{P} 0,$$

where $K$ is some universal constant. This together with (A.5) ensures that

$$(\boldsymbol{\alpha}' \mathbf{D}^{-1} \boldsymbol{\alpha})^{-1} \max \left| \frac{\alpha_j}{\sigma_j^2} (\hat{\epsilon}_{2j}) - \frac{\alpha_j}{\hat{\sigma}_j^2} (\hat{\epsilon}_{2j}) \right| = o_P(1).$$

Hence,

(A.6)
$$I_3 = \tilde{I}_3 + \boldsymbol{\alpha}' \mathbf{D}^{-1} \boldsymbol{\alpha} o_P(1).$$

Now we only need to consider $\tilde{I}_3$. Note that $\tilde{I}_3 = \sum \frac{\alpha_j}{\sigma_j^2} \hat{\epsilon}_{2j} \sim N(0, \frac{1}{n_2} \boldsymbol{\alpha}' \mathbf{D}^{-1} \times \boldsymbol{\Sigma} \mathbf{D}^{-1} \boldsymbol{\alpha})$. Since the variance term can be bounded as

$$\boldsymbol{\alpha}' \mathbf{D}^{-1} \boldsymbol{\Sigma} \mathbf{D}^{-1} \boldsymbol{\alpha} \leq \lambda_{\max}(\mathbf{R}) \boldsymbol{\alpha}' \mathbf{D}^{-1} \boldsymbol{\alpha},$$

by the assumption that $n \boldsymbol{\alpha}' \mathbf{D}^{-1} \boldsymbol{\alpha} \to \infty$ and $\lambda_{\max}(\mathbf{R})$ is bounded, we have $\tilde{I}_3 = \frac{1}{2} \boldsymbol{\alpha}' \mathbf{D}^{-1} \boldsymbol{\alpha} o_P(1)$. Combining this with (A.6) leads to

$$I_3 = \frac{1}{2} \boldsymbol{\alpha}' \mathbf{D}^{-1} \boldsymbol{\alpha} o_P(1).$$

We now examine $I_4$ and $I_5$. By the similar proof to (A.3) above we have

$$I_4 = p/n_1 + o_P(\sqrt{np/(n_1 n_2)}) \quad \text{and} \quad I_5 = p/n_2 + o_P(\sqrt{np/(n_1 n_2)}).$$

Thus the numerator can be written as

$$(\boldsymbol{\mu}_1 - \widehat{\boldsymbol{\mu}})' \widehat{\mathbf{D}}^{-1} (\widehat{\boldsymbol{\mu}}_1 - \widehat{\boldsymbol{\mu}}_2)$$
$$= (1 + o_P(1)) \frac{1}{2} \sum \alpha_j^2 / \sigma_j^2 - (p/n_1 - p/n_2)/2 + o_P(\sqrt{np/(n_1 n_2)})$$

and by (A.4)

$$\tilde{\Psi} = \frac{\sqrt{n_1 n_2/(pn)} \sum \alpha_j^2/\sigma_j^2(1 + o_P(1)) + \sqrt{p/(nn_1 n_2)}(n_1 - n_2)}{2\{1 + (n_1 n_2/(pn)) \sum \alpha_j^2/\sigma_j^2(1 + o_P(1))\}^{1/2}}.$$

Since $\frac{ax}{\sqrt{1 + a^2 x}}$ is an increasing function of $x$ and $\sum \frac{\alpha_j^2}{\sigma_j^2} \geq C_p$, in view of (A.1) and the definition of the parameter space $\Gamma$, we have

$$W(\hat{\delta}) = 1 - \Phi\left(\frac{[n_1 n_2/(pn)]^{1/2} C_p \{1 + o_P(1) + p(n_1 - n_2)/(n_1 n_2 C_p)\}}{2\sqrt{b_0}\{1 + n_1 n_2/(pn)C_p(1 + o_P(1))\}^{1/2}}\right).$$

If $p/(nC_p) \to 0$, then $W(\hat{\delta}) = 1 - \Phi(\frac{1}{2}[n_1 n_2/(pnb_0)]^{1/2} C_p\{1 + o_P(1)\})$. Furthermore, if $\{\frac{n_1 n_2}{pn}\}^{1/2} C_p \to C_0$ with $C_0$ some constant, then

$$W(\hat{\delta}) \xrightarrow{P} 1 - \Phi\left(\frac{C_0}{2\sqrt{b_0}}\right).$$

This completes the proof. □

PROOF OF THEOREM 2.  Suppose we have a new observation $\mathbf{X}$ from class $\mathcal{C}_1$. Then the posterior classification error of using $\hat{\delta}_{\mathbf{a}}(\cdot)$ is

$$W(\delta_{\mathbf{a}}, \boldsymbol{\theta}) = E^{\mathbf{a}}[P(\delta_{\mathbf{a}}(\mathbf{X}) < 0 | \mathbf{Y}_{ki}, i = 1, \ldots, n_k, \ k = 1, 2, \mathbf{a})]$$
$$= 1 - E^{\mathbf{a}}\Phi(\Psi_{\mathbf{a}} \operatorname{sign}(\mathbf{a}'\hat{\boldsymbol{\mu}}_1 - \mathbf{a}'\hat{\boldsymbol{\mu}}_2)),$$

where $\Psi_{\mathbf{a}} = \frac{\mathbf{a}'\boldsymbol{\mu}_1 - \mathbf{a}'\hat{\boldsymbol{\mu}}}{\sqrt{\mathbf{a}'\boldsymbol{\Sigma}\mathbf{a}}}$, $\Phi(\cdot)$ is the standard Gaussian distribution function, and $E^{\mathbf{a}}$ means expectation taken with respect to $\mathbf{a}$. We are going to show that

$$(A.7) \qquad\qquad\qquad \Psi_{\mathbf{a}} \xrightarrow{P} 0,$$

which together with the continuity of $\Phi(\cdot)$ and the dominated convergence theorem gives

$$\lim_p E^{\mathbf{a}}\Phi(\Psi_{\mathbf{a}} \operatorname{sign}(\mathbf{a}'\hat{\boldsymbol{\mu}}_1 - \mathbf{a}'\hat{\boldsymbol{\mu}}_2)) = 1/2.$$

Therefore, the posterior error $W(\hat{\delta}_{\mathbf{a}}, \boldsymbol{\theta})$ is no better than the random guessing.

Now, let us prove (A.7). Note that the random vector $\mathbf{a}$ can be written as

$$\mathbf{a} = \mathbf{Z}/\|\mathbf{Z}\|,$$

where $\mathbf{Z}$ is a $p$-dimensional standard Gaussian distributed random vector, independent of all the observations $\mathbf{Y}_{ki}$ and $\mathbf{X}$. Therefore,

$$(A.8) \ \Psi_{\mathbf{a}} = \frac{\mathbf{a}'\boldsymbol{\mu}_1 - \mathbf{a}'\hat{\boldsymbol{\mu}}}{\sqrt{\mathbf{a}'\boldsymbol{\Sigma}\mathbf{a}}} = \frac{\mathbf{Z}'\boldsymbol{\alpha}/\sqrt{p} - \sqrt{n/(n_1 n_2 p)}\mathbf{Z}'[\sqrt{(n_1 n_2)/n}(\hat{\boldsymbol{\epsilon}}_1 + \hat{\boldsymbol{\epsilon}}_2)]}{2\sqrt{\mathbf{Z}'\boldsymbol{\Sigma}\mathbf{Z}/p}},$$

where $\boldsymbol{\alpha} = \boldsymbol{\mu}_1 - \boldsymbol{\mu}_2$ and $\hat{\boldsymbol{\epsilon}}_k = \frac{1}{n_k}\sum_{i=1}^{n_k}\boldsymbol{\epsilon}_{ki}$, $k = 1, 2$. By the singular value decomposition we have

$$\boldsymbol{\Sigma} = \mathbf{Q}'\mathbf{V}\mathbf{Q},$$

where $\mathbf{Q}$ is an orthogonal matrix and $\mathbf{V} = \text{diag}\{\lambda_1, \ldots, \lambda_p\}$ is a diagonal matrix. Let $\tilde{\mathbf{Z}} = \mathbf{Q}\mathbf{Z}$, then $\tilde{\mathbf{Z}}$ is also a $p$-dimensional standard Gaussian random vector. Hence the denominator of $\Psi_{\mathbf{a}}$ can be written as

$$2\sqrt{\mathbf{Z}'\boldsymbol{\Sigma}\mathbf{Z}/p} = 2\left(\frac{1}{p}\sum_{j=1}^p \lambda_j \tilde{Z}_j^2\right)^{1/2},$$

where $\tilde{Z}_j$ is the $j$th entry of $\tilde{\mathbf{Z}}$. Since it is assumed that $\lim_p \frac{1}{p^2}\sum_{j=1}^p \lambda_j^2 < \infty$ and $\lim_p \frac{1}{p}\sum_{j=1}^p \lambda_j = \tau$ for some positive constant $\tau$, by the weak law of large numbers, we have

$$(A.9) \qquad \frac{1}{p}\sum_{j=1}^G \lambda_j \tilde{Z}_j^2 \xrightarrow{P} \tau.$$

Next, we study the numerator of $\Psi_{\mathbf{a}}$ in (A.8). Since $\frac{1}{p}\sum_{j=1}^p \alpha_j^2 \to 0$, the first term of the numerator converges to 0 in probability, that is,

$$(A.10) \qquad \mathbf{Z}'\boldsymbol{\alpha}/\sqrt{p} \xrightarrow{P} 0.$$

Let $\boldsymbol{\varepsilon} = \sqrt{\frac{n_1 n_2}{n}}(\hat{\boldsymbol{\epsilon}}_1 + \hat{\boldsymbol{\epsilon}}_2)$ and $\tilde{\boldsymbol{\varepsilon}} = \mathbf{V}^{-1/2}\mathbf{Q}\boldsymbol{\varepsilon}$, then $\tilde{\boldsymbol{\varepsilon}}$ has distribution $N(0, \mathbf{I})$ and is independent of $\tilde{\mathbf{Z}}$. The second term of the numerator can be written as

$$\mathbf{Z}'[\sqrt{n_1 n_2/n}(\hat{\boldsymbol{\epsilon}}_1 + \hat{\boldsymbol{\epsilon}}_2)] = \tilde{\mathbf{Z}}'\mathbf{V}^{1/2}\tilde{\boldsymbol{\varepsilon}} = \sum_{j=1}^p \sqrt{\lambda_j}\tilde{Z}_j \tilde{\varepsilon}_j.$$

Since $\frac{n}{n_1 n_2 p}\sum_{j=1}^p \lambda_j \to 0 < \infty$, it follows from the weak law of large number that

$$\sqrt{\frac{n}{n_1 n_2 p}}\sum_{j=1}^p \sqrt{\lambda_j}\tilde{Z}_j \tilde{\varepsilon}_j \xrightarrow{P} 0.$$

This together with (A.8), (A.9) and (A.10) completes the proof. $\square$

We need the following two lemmas to prove Theorem 3.

LEMMA A.1 [Cao (2007)]. *Let $n = n_1 + n_2$. Assume that there exist $0 < c_1 \le c_2 < 1$ such that $c_1 \le n_1/n_2 \le c_2$. Let $\tilde{T}_j = T_j - \frac{\mu_{j1} - \mu_{j2}}{\sqrt{S_{1j}^2/n_1 + S_{1j}^2/n_2}}$. Then for any $x \equiv x(n_1, n_2)$ satisfying $x \to \infty$ and $x = o(n^{1/2})$,*

$$\log P(\tilde{T}_j \ge x) \sim -x^2/2 \qquad \text{as } n_1, n_2 \to \infty.$$

*If in addition, if we have only $E|Y_{1ij}|^3 < \infty$ and $E|Y_{2ij}|^3 < \infty$, then*

$$\frac{P(\tilde{T}_j \geq x)}{1 - \Phi(x)} = 1 + O(1)(1+x)^3 n^{-1/2} d^3 \qquad \text{for } 0 \leq x \leq n^{1/6}/d,$$

*where $d = (E|Y_{1ij}|^3 + E|Y_{2ij}|^3)/(\mathrm{var}(Y_{1ij}) + \mathrm{var}(Y_{2ij}))^{3/2}$ and $O(1)$ is a finite constant depending only on $c_1$ and $c_2$. In particular,*

$$\frac{P(\tilde{T}_j \geq x)}{1 - \Phi(x)} \to 1$$

*uniformly in $x \in (0, o(n^{1/6}))$.*

LEMMA A.2. *Suppose Condition [1](b) holds and $\log p = o(n)$. Let $S_{kj}^2$ be the sample variance defined in Section [1](#), and $\sigma_{kj}^2$ be the variance of the $j$th feature in class $\mathcal{C}_k$. Suppose $\min \sigma_{kj}^2$ is bounded away from 0. Then we have the following uniform convergence result*

$$\max_{k=1,2, \; j=1,\ldots,p} |S_{kj}^2 - \sigma_{kj}^2| \overset{P}{\longrightarrow} 0.$$

PROOF. For any $\varepsilon > 0$, we know when $n_k$ is very large,

$$
\begin{aligned}
P&\bigg(\max_{k=1,2, \; j=1,\ldots,p} |S_{kj}^2 - \sigma_{kj}^2| > \varepsilon\bigg) \\
&\leq \sum_{k=1,2} \sum_{j=1}^{p} P(|S_{kj}^2 - \sigma_{kj}^2| > \varepsilon) \\
&\leq \sum_{k=1,2} \sum_{j=1}^{s} P\bigg(\bigg|\sum_{i=1}^{n_k} (\epsilon_{kij}^2 - \sigma_{kj}^2)\bigg| > n_k \varepsilon/3\bigg) \\
&\quad + \sum_{k=1,2} \sum_{j=1}^{s} P\bigg(\bigg|\sum_{i=1}^{n_k} \epsilon_{kij}\bigg| > n_k \sqrt{\varepsilon}/2\bigg) \\
&\equiv I_1 + I_2.
\end{aligned}
$$

$\text{(A.11)}$

It follows from Bernstein's inequality that

$$P\bigg(\bigg|\sum_{i=1}^{n_k}(\epsilon_{kij}^2 - \sigma_{kj}^2)\bigg| > n_k \varepsilon/3\bigg) \leq 2\exp\bigg\{-\frac{1}{2}\frac{n_k^2 \varepsilon^2}{9\nu_1 + 3M_1 n_k \varepsilon}\bigg\}$$

and

$$P\bigg(\bigg|\sum_{i=1}^{n_k} \epsilon_{kij}\bigg| > n_k \sqrt{\varepsilon}/2\bigg) \leq 2\exp\bigg\{-\frac{1}{2}\frac{n_k^2 \varepsilon}{4\nu_2 + 2M_2 n_k \sqrt{\varepsilon}}\bigg\}.$$

Since $\log p = o(n)$, we have $I_1 = o(1)$ and $I_2 = o(1)$. These together with (A.11) completes the proof of Lemma A.2. $\quad\square$

PROOF OF THEOREM 3. We divide the proof into two parts. (a) Let us first look at the probability $P(\max_{j>s}|T_j| > x)$. Clearly,

$$(A.12) \qquad P\left(\max_{j>s}|T_j| > x\right) \le \sum_{j=s+1}^{p} P(|T_j| \ge x).$$

Note that for all $j > s$, $\alpha_j = \mu_{j1} - \mu_{j2} = 0$. By Condition 1(b) and Lemma A.1, the following inequality holds for $0 \le x \le n^{1/6}/d$,

$$P(T_j \ge x) = (1 - \Phi(x))(1 + C(1+x)^3 n^{-1/2} d^3),$$

where $C$ is a constant that only depends on $c_1$ and $c_2$, and

$$d = (E|Y_{1ij}|^3 + E|Y_{2ij}|^3)/(\sigma_{1j}^2 + \sigma_{2j}^2)^{3/2}$$

with $\sigma_{kj}^2$ the $j$th diagonal element of $\Sigma_k$. For the normal distribution, we have the following tail probability inequality

$$1 - \Phi(x) \le \frac{1}{\sqrt{2\pi}}\frac{1}{x}e^{-x^2/2}.$$

This together with the symmetry of $T_j$ gives

$$P(|T_j| \ge x) \le 2\frac{1}{\sqrt{2\pi}}\frac{1}{x}e^{-x^2/2}(1 + C(1+x)^3 n^{-1/2} d^3).$$

Combining the above inequality with (A.12), we have

$$\sum_{j>s} P(|T_j| \ge x) \le (p-s)\frac{2}{\sqrt{2\pi}}\frac{1}{x}e^{-x^2/2}(1 + C(1+x)^3 n^{-1/2} d^3).$$

Since $\log(p-s) = o(n^\gamma)$ with $0 < \gamma < \frac{1}{3}$, if we let $x \sim cn^{\gamma/2}$, then

$$\sum_{j>s} P(|T_j| \ge x) \to 0,$$

which along with (A.12) yields

$$P\left(\max_{j>s}|T_j| > x\right) \to 0.$$

(b) Next, we consider $P(\min_{j\le s}|T_j| \le x)$. Notice that for $j \le s$, $\alpha_j = \mu_{1j} - \mu_{2j} \ne 0$. Let $\eta_j = \frac{\alpha_j}{\sqrt{S_{1j}^2/n_1 + S_{1j}^2/n_2}}$ and define

$$\tilde{T}_j = T_j - \eta_j.$$

Then following the same lines as those in (a), we have

$$\sum_{j \leq s} P(|\tilde{T}_j| \geq x) \leq s \frac{2}{\sqrt{2\pi}} \frac{1}{x} e^{-x^2/2}(1 + C(1+x)^3 n^{-1/2} d^3) \to 0.$$

It follows from Lemma A.2 that,

$$\max_{j \leq s} |S_{kj}^2 - \sigma_{kj}^2| \xrightarrow{P} 0, \qquad k = 1, 2.$$

Hence, uniformly over $j = 1, \ldots, s$, we have

$$\eta_j = \frac{|\alpha_j|}{\sqrt{\sigma_{1j}^2/n_1 + \sigma_{2j}^2/n_2}}(1 + o_P(1)).$$

Therefore,

$$\min_{j \leq s} |\eta_j| = \min_{j \leq s} \frac{\sqrt{n_1}|\alpha_j|}{\sqrt{\sigma_{1j}^2 + \sigma_{2j}^2 n_1/n_2}}(1 + o_P(1)) \geq \min_{j \leq s} \frac{\sqrt{n_1}|\alpha_j|}{\sqrt{\sigma_{1j}^2 + c_2 \sigma_{2j}^2}}(1 + o_P(1))$$

with $c_2$ defined in Theorem 3. Let $\alpha_0 = \min_{j \leq s} |\mu_{j1} - \mu_{j2}|/\sqrt{\sigma_{1j}^2 + c_2\sigma_{2j}^2}$. Then it follows that

$$P\left(\min_{j \leq s} |T_j| \leq x\right) \leq P\left(\max_{j \leq s} |\tilde{T}_j| \geq \min_{j \leq s} |\eta_j| - x\right)$$

$$\leq P\left(\max_{j \leq s} |\tilde{T}_j| \geq \sqrt{n_1}\alpha_0(1 + o_P(1)) - x\right).$$

By part (a), we know that $x \sim cn^{\gamma/2}$ and $\log(p - s) = o(n^\gamma)$. Thus if $\alpha_0 \sim \min_{j \leq s} \frac{|\mu_{j1} - \mu_{j2}|}{\sqrt{\sigma_{j1}^2 + \sigma_{j2}^2}} = n^{-\gamma}\beta_n$ for some $\beta_n \to \infty$, then similarly to part (a), we have

$$P\left(\min_{j \leq s} |T_j| \leq x\right) \to 0.$$

Combination of part (a) and part (b) completes the proof. □

PROOF OF THEOREM 4. The classification error of the truncated classifier $\hat{\delta}_{\mathrm{NC}}^m$ is

$$W(\hat{\delta}_{\mathrm{NC}}^m, \boldsymbol{\theta}) = 1 - \Phi\left(\frac{\sum_{j=1}^m \hat{\alpha}_j(\mu_{1j} - \hat{\mu}_j)}{\sum_{j=1}^m \hat{\alpha}_j^2}\right).$$

We first consider the denominator. Note that $\hat{\alpha}_j \sim N(\alpha_j, \frac{n}{n_1 n_2})$. It can be shown that

$$\left(\frac{4n}{n_1 n_2}\sum_{j=1}^m \alpha_j^2 + \frac{2mn^2}{n_1^2 n_2^2}\right)^{-1/2}\sum_{j=1}^m\left(\hat{\alpha}_j^2 - \alpha_j^2 - \frac{n}{n_1 n_2}\right) \xrightarrow{D} N(0, 1),$$

which together with the assumption $\frac{n}{\sqrt{m}} \sum_{j=1}^{m} \alpha_j^2 \to \infty$ gives

$$\sum_{j=1}^{m} \hat{\alpha}_j^2 = \sum_{j=1}^{m} \alpha_j^2 + \frac{mn}{n_1 n_2} + \left\{ \frac{4n}{n_1 n_2} \sum_{j=1}^{m} \alpha_j^2 + \frac{2mn^2}{n_1^2 n_2^2} \right\}^{1/2} O_p(1)$$

$$= (1 + o_P(1)) \sum_{j=1}^{m} \alpha_j^2 + \frac{mn}{n_1 n_2}.$$

Next, let us look at the numerator. We decompose it as

(A.13) $$\sum_{j=1}^{m} \hat{\alpha}_j (\mu_{1j} - \hat{\mu}_j) = \frac{1}{2} \sum_{j=1}^{m} \alpha_j^2 - \sum_{j=1}^{m} \alpha_j \hat{\epsilon}_{2j} - \frac{1}{2} \sum_{j=1}^{m} (\hat{\epsilon}_{1j}^2 - \hat{\epsilon}_{2j}^2).$$

Since the second term above has the distribution $N(0, \sum_{j=1}^{m} \alpha_j^2 / n_2)$, it follows from the assumption $n \sum_{j=1}^{m} \alpha_j^2 \to \infty$ that

$$\sum_{j=1}^{m} \alpha_j \hat{\epsilon}_{2j} = o_P(1) \sum_{j=1}^{m} \alpha_j^2.$$

The third term in (A.13) can be written as

$$\sum_{j=1}^{m} (\hat{\epsilon}_{1j}^2 - \hat{\epsilon}_{2j}^2) = \frac{m}{n_1} - \frac{m}{n_2} + O_p\left( \frac{nm}{n_1 n_2} \right) = \frac{m(n_2 - n_1)}{n_1 n_2} + o_P(1) \sum_{j=1}^{m} \alpha_j^2.$$

Hence the numerator is

$$\sum_{j=1}^{m} \hat{\alpha}_j (\mu_{1j} - \hat{\mu}_j) = \frac{m(n_2 - n_1)}{n_1 n_2} + (1 + o_P(1)) \sum_{j=1}^{m} \alpha_j^2.$$

Therefore, the classification error is

$$W(\hat{\delta}_{\mathrm{NC}}^m, \boldsymbol{\theta}) = 1 - \Phi\left( \frac{(1 + o_P(1)) \sum_{j=1}^{m} \alpha_j^2 + m(n_1 - n_2)/(n_1 n_2)}{2\{(1 + o_P(1)) \sum_{j=1}^{m} \alpha_j^2 + mn/(n_1 n_2)\}^{1/2}} \right).$$

This concludes the proof. $\square$

PROOF OF THEOREM 5. Note that the classification error of $\hat{\delta}_{\mathrm{FAIR}}^{b_n}$ is

$$W(\hat{\delta}_{\mathrm{FAIR}}^{b_n}(\mathbf{x}), \boldsymbol{\theta}) = 1 - \Phi\left( \frac{\sum_j (\mu_{1j} - \hat{\mu}_j) \hat{\alpha}_j 1\{|\hat{\alpha}_j| \geq b_n\}}{\sum_j \hat{\alpha}_j^2 1\{|\hat{\alpha}_j| \geq b_n\}} \right) \equiv 1 - \Phi(\Psi^H).$$

We divide the proof into two parts: the numerator and the denominator.
(a) First, we study the numerator of $\Psi^H$. It can be decomposed as

$$\sum_j (\mu_{1j} - \hat{\mu}_j) \hat{\alpha}_j 1\{|\hat{\alpha}_j| \geq b_n\} = I_1 + I_2,$$

where $I_1 = \sum_{j \in \mathcal{A}} (\mu_{1j} - \hat{\mu}_j) \hat{\alpha}_j 1\{|\hat{\alpha}_j| \geq b_n\}$ and $I_2 = \sum_{j \in \mathcal{A}^c} (\mu_{1j} - \hat{\mu}_j) \hat{\alpha}_j \times 1\{|\hat{\alpha}_j| \geq b_n\}$ with $\mathcal{A}^c$ the complementary of the set $\mathcal{A}$. Note that

$$
\begin{aligned}
I_2 &= \tfrac{1}{2} \sum_{j \in \mathcal{A}^c} \alpha_j^2 1\{|\hat{\alpha}_j| \geq b_n\} - \sum_{j \in A^c} \alpha_j \hat{\epsilon}_{2j} 1\{|\hat{\alpha}_j| \geq b_n\} \\
&\quad - \tfrac{1}{2} \sum_{j \in A^c} (\hat{\epsilon}_{1j}^2 - \hat{\epsilon}_{2j}^2) 1\{|\hat{\alpha}_j| \geq b_n\} \\
&\equiv \tfrac{1}{2} I_{2,1} - I_{2,2} - \tfrac{1}{2} I_{2,3}.
\end{aligned}
$$

Since $\hat{\alpha}_j \sim N(\alpha_j, \frac{n}{n_1 n_2})$, it follows from the normal tail probability inequality that for every $j \in \mathcal{A}^c$ and $b_n > \max_{j \in \mathcal{A}^c} |\alpha_j|$,

$$
\begin{aligned}
\text{(A.14)} \qquad P(|\hat{\alpha}_j| \geq b_n) &\leq P\left(|\hat{\alpha}_j - \alpha_j| \geq b_n - \max_{j \in \mathcal{A}^c} |\alpha_j|\right) \\
&\leq M \frac{\exp\{-n_1 n_2 (b_n - \max_{j \in A^c} |\alpha_j|)^2 / (2n)\}}{\sqrt{n_1 n_2 n^{-1}} (b_n - \max_{j \in \mathcal{A}^c} |\alpha_j|)},
\end{aligned}
$$

where $M$ is a generic constant. Thus for every $\varepsilon > 0$, if $\log(p - m)/[n(b_n - \max_{j \in \mathcal{A}^c} |\alpha_j|)^2] \to 0$ and $\max_{j \in \mathcal{A}^c} |\alpha_j| < b_n$, we have

$$
\begin{aligned}
P(|I_{2,1}| \geq \varepsilon) &\leq \varepsilon^{-1} \sum_{j \in A^c} \alpha_j^2 P(|\hat{\alpha}_j| \geq b_n) \\
&\leq M \max_{j \in \mathcal{A}^c} \alpha_j^2 \frac{(p - m)}{\varepsilon} \frac{\exp\{-n_1 n_2 (b_n - \max_{j \in A^c} |\alpha_j|)^2 / (2n)\}}{\sqrt{n_1 n_2 n^{-1}} (b_n - \max_{j \in \mathcal{A}^c} |\alpha_j|)},
\end{aligned}
$$

which tends to zero. Hence,

$$
\text{(A.15)} \qquad\qquad\qquad I_{2,1} \xrightarrow{P} 0.
$$

We next consider $I_{2,2}$. Since $E(\hat{\epsilon}_{2j})^2 = \frac{1}{n_2}$, $\log(p - m)/[n(b_n - \max_{j \in A^c} |\alpha_j|)^2] \to 0$, and $\max_{j \in \mathcal{A}^c} |\alpha_j| < b_n$, we have

$$
\begin{aligned}
P(|I_{2,2}| \geq \varepsilon) &\leq \varepsilon^{-1} \sum_{j \in A^c} E|\hat{\epsilon}_{2j} \alpha_j 1\{|\hat{\alpha}_j| \geq b_n\}| \\
&\leq \varepsilon^{-1} \sum_{j \in A^c} \{E(\hat{\epsilon}_{2j})^2\}^{1/2} \{E|\alpha_j^2 1\{|\hat{\alpha}_j| \geq b_n\}|\}^{1/2} \\
&\leq M \frac{(p - m) \max_{j \in \mathcal{A}^c} |\alpha_j|}{\varepsilon \sqrt{n_2}} \frac{\exp\{-n_1 n_2 (b_n - \max_{j \in A^c} |\alpha_j|)^2 / (4n)\}}{\sqrt{n_1 n_2 n^{-1}} (b_n - \max_{j \in \mathcal{A}^c} |\alpha_j|)},
\end{aligned}
$$

which converges to 0. Therefore,

$$
\text{(A.16)} \qquad\qquad\qquad I_{2,2} \xrightarrow{P} 0.
$$

Then, we consider $I_{2,3}$. Since $c_1 \le n_1/n_2 \le c_2$ and $E(\hat{\epsilon}_{1j}^2 - \hat{\epsilon}_{2j}^2)^2 = \frac{3n_1^2 + 3n_2^2 - 2n_1 n_2}{n_1^2 n_2^2} \le \frac{3c_2 + 3 - 2c_1}{c_1 n_2^2}$, by (A.14) we have for every $\varepsilon > 0$,

$$P(|I_{2,3}| \ge \varepsilon) \le \varepsilon^{-1} E \left| \sum_{j \in A^c} (\hat{\epsilon}_{1j}^2 - \hat{\epsilon}_{2j}^2) 1\{|\hat{\alpha}_j| \ge b_n\} \right|$$

$$\le \varepsilon^{-1} \sum_{j \in A^c} \{E(\hat{\epsilon}_{1j}^2 - \hat{\epsilon}_{2j}^2)^2\}^{1/2} P(|\hat{\alpha}_j| \ge b_n)^{1/2}$$

$$\le M \frac{\sum_{j \in A^c} P(|\hat{\alpha}_j| \ge b_n)}{n_2 \varepsilon} \to 0,$$

where $M$ is some generic constant. Thus, $I_{2,3} \xrightarrow{P} 0$. Combination of this with (A.15) and (A.16) entails

$$I_2 = o_P(1).$$

We now deal with $I_1$. Decompose $I_1$ similarly as

$$I_1 = \sum_{j \in A} (\mu_{1j} - \hat{\mu}_{1j}) \hat{\alpha}_j 1\{|\hat{\alpha}_j| \ge b_n\} + \frac{1}{2} \sum_{j \in A} \hat{\alpha}_j^2 1\{|\hat{\alpha}_j| \ge b_n\}$$

$$\ge \sum_{j \in A} (\mu_{1j} - \hat{\mu}_{1j}) \hat{\alpha}_j 1\{|\hat{\alpha}_j| \ge b_n\} + \frac{1}{2} \sum_{j \in A} (\hat{\alpha}_j^2 - b_n^2)$$

$$\equiv I_{1,1} + \tfrac{1}{2} I_{1,2}.$$

We first study $I_{1,2}$. By using $\hat{\alpha}_j \sim N(\alpha_j, \frac{n}{n_1 n_2})$, it can be shown that

$$(A.17) \quad \left( \frac{4n}{n_1 n_2} \sum_{j \in A} \alpha_j^2 + \frac{2mn^2}{n_1^2 n_2^2} \right)^{-1/2} \sum_{j \in A} \left( \hat{\alpha}_j^2 - \alpha_j^2 - \frac{n}{n_1 n_2} \right) \xrightarrow{D} N(0,1).$$

Since $\frac{n}{\sqrt{m}} \sum_{j=1}^m \alpha_j^2 \to \infty$, we have $(\frac{4n}{n_1 n_2} \sum_{j \in A} \alpha_j^2 + \frac{2mn^2}{n_1^2 n_2^2})^{1/2} / \sum_{j \in A} \alpha_j^2 \to 0$. Therefore,

$$I_{1,2} = \sum_{j \in A} (\alpha_j^2 - b_n^2) + \frac{nm}{n_1 n_2} + \left( \frac{4n}{n_1 n_2} \sum_{j \in A} \alpha_j^2 + \frac{2mn^2}{n_1^2 n_2^2} \right)^{1/2} O_p(1)$$

$$= (1 + o_P(1)) \sum_{j \in A} \alpha_j^2 + \frac{nm}{n_1 n_2} - mb^2.$$

Next, we look at $I_{1,1}$. For any $\varepsilon > 0$,

$$P(|I_{1,1}| \ge \varepsilon) \le \frac{1}{\varepsilon} E|I_{1,1}| \le \frac{1}{\varepsilon} \sum_{j \in A} \{E|\mu_{1j} - \hat{\mu}_{1j}|^2 E|\hat{\alpha}_j|^2\}^{1/2}$$

$$= \frac{1}{\sqrt{n_1} \varepsilon} \sum_{j \in A} \sqrt{\alpha_j^2 + n/(n_1 n_2)}.$$

When $n$ is large enough, the above probability can be bounded by

$$P(|I_{1,1}| \geq \varepsilon) \leq \sqrt{2/(n_1 \varepsilon^2)} \sum_{j \in \mathcal{A}} |\alpha_j|,$$

which along with the assumption $\sum_{j \in \mathcal{A}} |\alpha_j| / [\sqrt{n} \sum_{j \in \mathcal{A}} \alpha_j^2] \to 0$ gives

$$I_{1,1} = o_P(1) \sum_{j \in \mathcal{A}} \alpha_j^2.$$

It follows that the numerator is bounded from below by

$$(1 + o_P(1)) \frac{1}{2} \sum_{j \in \mathcal{A}} \alpha_j^2 + \frac{mn}{2n_1 n_2} - \frac{1}{2} mb^2.$$

(b) Now, we study the denominator of $\Psi$. Let

$$\sum_j \hat\alpha_j^2 1\{|\hat\alpha_j| \geq b_n\} = \sum_{j \in \mathcal{A}} \hat\alpha_j^2 1\{|\hat\alpha_j| \geq b_n\} + \sum_{j \in \mathcal{A}^c} \hat\alpha_j^2 1\{|\hat\alpha_j| \geq b_n\} \equiv J_1 + J_2.$$

We first show that $J_2 \xrightarrow{P} 0$. Note that $E\hat\alpha_j^4 = \alpha_j^4 + 6n(n_1 n_2)^{-1} \alpha_j^2 + 3n^2(n_1 n_2)^{-2}$. Thus,

$$P(|J_2| \geq \varepsilon) \leq \frac{1}{\varepsilon} E|J_2| = \sum_{j \in A^c} E\hat\alpha_j^2 1\{|\hat\alpha_j| \geq b_n\} / \varepsilon \leq \frac{1}{\varepsilon} \sum_{j \in A^c} \{E\hat\alpha_j^4 P(|\hat\alpha_j| \geq b_n)\}^{1/2}$$

$$\leq \frac{1}{\varepsilon} \sum_{j \in A^c} \{(\alpha_j^4 + 6n(n_1 n_2)^{-1} \alpha_j^2 + 3n^2(n_1 n_2)^{-2}) P(|\hat\alpha_j| \geq b_n)\}^{1/2}.$$

This together with (A.14) and the assumption that $\log(p - m)/[n(b_n - \max_{j \in \mathcal{A}^c} |\alpha_j|)^2] \to 0$ yields $J_2 \xrightarrow{P} 0$ as $n \to \infty$, $p \to \infty$. Now we study term $J_1$. By (A.17), we have

$$J_1 \leq \sum_{j \in \mathcal{A}} \hat\alpha_j^2 = (1 + o_P(1)) \sum_{j \in \mathcal{A}} \alpha_j^2 + \frac{mn}{n_1 n_2}.$$

Hence the denominator is bounded from above by $(1 + o_P(1)) \sum_{j \in \mathcal{A}} \alpha_j^2 + \frac{mn}{n_1 n_2}$. Therefore,

$$\Psi^H \geq \frac{(1 + o_P(1)) \sum_{j \in \mathcal{A}} \alpha_j^2 + (mn/(n_1 n_2)) - mb^2}{2\sqrt{(1 + o_P(1)) \sum_{j \in \mathcal{A}} \alpha_j^2 + (mn/(n_1 n_2))}}.$$

It follows that the classification error is bounded from above by

$$1 - \Phi\left( \frac{(1 + o_P(1)) \sum_{j \in \mathcal{A}} \alpha_j^2 + (mn/(n_1 n_2)) - mb^2}{2\sqrt{(1 + o_P(1)) \sum_{j \in \mathcal{A}} \alpha_j^2 + (mn/(n_1 n_2))}} \right).$$

This completes the proof. $\quad\square$

**Acknowledgments.**   The authors acknowledge gratefully the helpful comments of referees that led to the improvement of the presentation and the results of the paper.

DEPARTMENT OF OPERATIONS RESEARCH
 AND FINANCIAL ENGINEERING
PRINCETON UNIVERSITY
PRINCETON, NEW JERSEY 08544
USA
E-MAIL: jqfan@princeton.edu

INFORMATION AND OPERATIONS
 MANAGEMENT DEPARTMENT
MARSHALL SCHOOL OF BUSINESS
UNIVERSITY OF SOUTHERN CALIFORNIA
LOS ANGELES, CALIFORNIA 90089
USA
E-MAIL: fanyingy@marshall.usc.edu